\documentclass[11pt]{article}

\usepackage{amsmath,amssymb,amsthm,enumitem}
\usepackage{pstricks,pst-node,xy,float}
\xyoption{all}

\numberwithin{equation}{section}
\newtheorem{thm}{Theorem}[section] 
\newtheorem{prp}[thm]{Proposition}
\newtheorem{lmm}[thm]{Lemma}   
\newtheorem{crl}[thm]{Corollary} 
\newtheorem{dfn}[thm]{Definition} 
\theoremstyle{definition}

\newtheorem{rmk}[thm]{Remark}

\def\ov#1{\overline{#1}}
\def\un#1{\underline{#1}}
\def\ti#1{\tilde{#1}}

\def\sf#1{\textsf{#1}}
\def\sm#1{\begin{small}{#1}\end{small}}
\def\e_ref#1{(\ref{#1})}
\def\lr#1{\langle{#1}\rangle}
\def\llrr#1{\langle\!\langle{#1}\rangle\!\rangle}
\def\blr#1{\big\langle{#1}\big\rangle}

\def\bblr#1{\bigg\langle{#1}\bigg\rangle}
\def\BE#1{\begin{equation}\label{#1}}
\def\EE{\end{equation}}

\def\lra{\longrightarrow}
\def\Lra{\Longrightarrow}
\def\Llra{\Longleftrightarrow}

\def\C{\mathbb C}
\def\cD{\mathcal D}
\def\d{\textnormal d}

\def\nh{\textnormal h}
\def\fj{\mathfrak j}
\def\I{\mathbb I}
\def\fI{\mathfrak i}

\def\N{\mathcal N}

\def\R{\mathbb R}
\def\cR{\mathcal R}

\def\al{\alpha}
\def\be{\beta}
\def\de{\delta}
\def\ep{\epsilon}
\def\ga{\gamma}
\def\io{\iota}
\def\na{\nabla}
\def\om{\omega}
\def\si{\sigma}
\def\th{\theta}

\def\vph{\varphi}

\def\ze{\zeta}

\def\Ga{\Gamma}

\def\Si{\Sigma}
\def\Th{\Theta}

\def\d{\textnormal{d}}
\def\E{\textnormal{e}}
\def\diam{\textnormal{diam}}
\def\dim{\textnormal{dim}}

\def\End{\textnormal{End}}

\def\Hom{\textnormal{Hom}}
\def\id{\textnormal{id}}
\def\Id{\textnormal{Id}}
\def\Im{\textnormal{Im}}

\def\Mat{\textnormal{Mat}}
\def\rk{\textnormal{rk}}

\def\SO{\textnormal{SO}}

\def\supp{\textnormal{supp}}

\def\nv{\textnormal{v}}

\def\U{\textnormal{U}}

\def\dbar{\bar\partial}

\def\i{\infty}

\begin{document}

\title{Basic Riemannian Geometry and Sobolev Estimates\\
used in Symplectic Topology}
\author{Aleksey Zinger\thanks{Partially supported by DMS grant 0846978}}
\maketitle

\begin{abstract}
This note collects a number of standard statements in Riemannian geometry
and in Sobolev-space theory
that play a prominent role in analytic approaches to symplectic topology.
These include relations between connections and complex structures,
estimates on exponential-like maps, and dependence of constants
in Sobolev and elliptic estimates.
\end{abstract}

\tableofcontents

\section{Connections in real vector bundles}

\subsection{Connections and splittings}
\label{genconn_subs}

\noindent
Suppose $M$ is a smooth manifold and $\pi\!:E\!\lra\!M$ is a vector bundle.
Trivializations of $M$ induce a bundle inclusion $\pi^*E\!\lra\!TE$
so that the sequence of vector bundles over~$E$ 
\begin{equation}\label{bundleseq_e1}
0 \lra \pi^*E \lra TE \stackrel{\d\pi}{\lra} \pi^*TM\lra 0
\end{equation}
is exact.
For each $f\!\in\!C^{\i}(M)$, define 
\begin{equation}\label{mf_dfn_e}
m_f\!: E\lra E \qquad\hbox{by}\qquad
m_f(v)=f\big(\pi(v))\cdot v \quad \forall~v\!\in\!E.
\end{equation}
We then have a commutative diagram
\begin{equation}\label{bundleseq_e2}
\xymatrix{0 \ar[r]& \pi^*E \ar[r] \ar[d]^{\pi^*m_f} & TE \ar[r]^{\d\pi} \ar[d]^{\d m_f} &
\pi^*TM \ar[r] \ar[d]^{\id} & 0\\
0 \ar[r]& \pi^*E \ar[r]& m_f^*TE \ar[r]^{\d\pi} & \pi^*TM \ar[r]  & 0}
\end{equation}
of bundle maps over $E$.\\

\noindent
A \sf{connection~in~$E$} is an $\R$-linear map
\begin{equation}\label{na_dfn_e}
\na\!:\Ga(M;E)\lra\Ga(M;T^*M\!\otimes_{\R}\!E) \quad\hbox{s.t.}\quad
\na(f\xi)=\d f\!\otimes\!\xi+f\na\xi 
~~~\forall~f\!\in\!C^{\i}(M),~\xi\!\in\!\Ga(M;E).
\end{equation}
The Leibnitz property implies that any two connections in $E$ differ
by a one-form on~$M$.
In other words, if $\na$ and $\ti\na$ are connections in $E$ 
there exists
\begin{gather}
\th\in\Ga\big(M;T^*M\!\otimes_{\R}\!\Hom_{\R}(E,E)\big) \qquad\hbox{s.t.}\notag\\
\label{connection_diff_e}
\ti\na_v\xi=\na_v\xi+\big\{\th(v)\big\}\xi \qquad\forall~
\xi\!\in\!\Ga(M;E),~v\!\in\!T_xM,~x\!\in\!M.
\end{gather}\\

\noindent
A connection $\na$ in $E$ is  a local operator, 
i.e.~the value of $\na\xi$ at a point $x\!\in\!M$ depends only on the restriction 
of $\xi$ to any neighborhood $\U$ of~$x$.
If $f$ is a smooth function on $M$ supported in $\U$ such that $f(x)\!=\!1$, then
\begin{equation}\label{local_prp_e}
\na\xi\big|_x=\na\big(f\xi)\big|_x-\d_xf\!\otimes\!\xi(x)
\end{equation}
by~\e_ref{na_dfn_e}. 
The right-hand side of~\e_ref{local_prp_e} depends only on $\xi|_{\U}$.\\

\noindent
In fact, a connection $\na$ in $E$ is a first-order differential operator.
Suppose $\U$ is an open subset of $M$ and $\xi_1,\ldots,\xi_n\!\in\!\Ga(\U;E)$ is
a frame for $E$ on~$\U$, i.e.
$$\xi_1(x),\ldots,\xi_n(x)\in E_x$$
is a basis for $E_x$ for all $x\!\in\!\U$.
By definition of $\na$, there exist 
$$\th^k_l\in \Ga(M;T^*M) \qquad\hbox{s.t.}\qquad
\na\xi_l=\sum_{k=1}^{k=n}\xi_k\th^k_l 
\equiv \sum_{k=1}^{k=n}\th^k_l\!\otimes\!\xi_k
\quad\forall~l\!=\!1,\ldots,n.$$
We will call
$$\th \equiv \big(\th^k_l\big)_{k,l=1,\ldots,n}
\in \Ga\big(\Si;T^*M\!\otimes_{\R}\!\Mat_n\R\big)$$
\sf{the connection one-form of $\na$ with respect to the frame $(\xi_k)_k$}.
For an arbitrary section
$$\xi=\sum_{l=1}^{l=n}f^l\xi_l\in \Ga(\U;E),$$
by~\e_ref{na_dfn_e} we have
\begin{gather}\label{na_e3}
\na\xi=\sum_{k=1}^{k=n}\xi_k \Big(\d f^k +
\sum_{l=1}^{l=n}\th^k_lf^l\Big), \qquad\hbox{i.e.}\qquad
\na\big(\un{\xi}\cdot\un{f}^t\big)=\un\xi\cdot\big\{\d+\th\big\}\un{f}^t,\\
\hbox{where}\qquad \un\xi=(\xi_1,\ldots,\xi_n), \quad \un{f}=(f^1,\ldots,f^n).
\label{Cna_e3b}
\end{gather}
Thus, $\na$ is a first-order differential operator.
It is immediate from~\e_ref{na_dfn_e} that the symbol of~$\na$ is given~by
$$\si_{\na}\!: T^*M\lra  \Hom\big(E,T^*M\!\otimes_{\R}\!E\big), \qquad
\big\{\si_{\na}(\eta)\big\}(f)=\eta\otimes f.$$\\

\noindent
Since $M\!\subset\!E$ as the zero section, there is a natural splitting
\begin{equation}\label{splitting_e0}
TE|_M\approx TM\oplus E
\end{equation}
of the exact sequence~\e_ref{bundleseq_e1} restricted to $M$.
If $x\!\in\!M$ and $\xi\!\in\!\Ga(M;E)$ is such that $\xi(x)\!=\!0$,
then
\begin{equation}\label{splitting_e1}
\na\xi\big|_x=\pi_2|_x\circ\d_x\xi\,,
\end{equation}
where $\pi_2|_x\!:T_xE\!\lra\!E_x$ is the projection onto the second 
component in~\e_ref{splitting_e0}.
This observation follows from~\e_ref{connection_diff_e},
as well as from~\e_ref{na_e3}.

\begin{lmm}\label{genconn_lmm}
Suppose $M$ is a smooth manifold and $\pi\!:E\!\lra\!M$ is a vector bundle.
A connection~$\na$ in $E$ induces a splitting 
\begin{equation}\label{genconn_e1}
TE\approx \pi^*TM\oplus\pi^*E
\end{equation}
of the exact sequence~\e_ref{bundleseq_e1} extending 
the splitting~\e_ref{splitting_e0} such that 
\begin{equation}\label{genconn_e2a}
\na\xi\big|_x=\pi_2|_x\circ\d_x\xi 
\qquad\forall~\xi\!\in\!\Ga(M;E),~x\!\in\!M,
\end{equation}
where $\pi_2|_x\!:T_xE\!\lra\!E_x$ is the projection onto the second 
component in~\e_ref{genconn_e1}, and
\begin{equation}\label{genconn_e2b}
\d m_t\approx \pi^*\id\oplus \pi^*m_t \qquad\forall~t\!\in\!\R,
\end{equation}
i.e.~the splitting is consistent with the commutative diagram~\e_ref{bundleseq_e2}.
\end{lmm}

\noindent
{\it Proof:} For each $x\!\in\!M$ and $v\!\in\!E_x$, choose $\xi\!\in\!\Ga(M;E)$
such that $\xi(x)\!=\!v$ and let 
$$T_vE^{\nh}=\Im\,\{\d\xi\!-\!\na\xi\}\big|_x\subset T_vE.$$
Since $\pi\!\circ\!\xi\!=\!\id_M$, 
$$\d_v\pi\circ\big\{\d\xi\!-\!\na\xi\big\}=\id_{T_xM}
\qquad\Lra\qquad
T_vE\approx T_vE^{\nh}\oplus E_x \approx T_xM\oplus E_x.$$
If $v\!=\!0$, then by~\e_ref{splitting_e1}
$$T_vE^{\nh}=T_vM.$$
If $v\!\neq\!0$, $\ze\!\in\!\Ga(M;E)$ is another section such that $\ze(x)\!=\!v$,
and $\U$ is sufficiently small,
then $\ze\!=\!f\xi$ for some $f\!\in\!C^{\i}(\U)$ with $f(x)\!=\!1$ and~thus
\begin{equation*}\begin{split}
\{\d\ze\!-\!\na\ze\big\}\big|_x &=\{\d(f\xi)\!-\!\na(f\xi)\big\}\big|_x
=\big\{\d_xf\!\otimes\!\xi(x)\!+\!f(x)\d_x\xi\big\}
-\big\{\d_xf\!\otimes\!\xi(x)\!+\!f(x)\na\xi|_x\big\}\\
&=\big\{\d\xi\!-\!\na\xi\big\}\big|_x\,.
\end{split}\end{equation*}
The second equality above is obtained by considering a trivialization of $E$ near $x$.
Thus, $T_vE^{\nh}$ is independent of the choice of $\xi$
and we obtain a well-defined splitting~\e_ref{genconn_e1} of~\e_ref{bundleseq_e1}
that satisfies~\e_ref{genconn_e2a} and extends~\e_ref{splitting_e0}.\\

\noindent
It remains to verify~\e_ref{genconn_e2b}. Since $\pi\!\circ\!m_t\!=\!\pi$,
$\d\pi\!\circ\!\d m_t\!=\!\d\pi$, i.e.~the first component of $\d m_t$ vanishes
on $TE$ and is the identity on~$\pi^*TM$.
On the other hand, if $\xi\!\in\!\Ga(M;E)$ and $x\!\in\!M$, then
\begin{equation*}\begin{split}
T_{t\xi(x)}E^{\nh} &\equiv 
\Im\big\{\d(m_t\!\circ\!\xi)\!-\!\na(t\xi)\big\}\big|_x
=\Im\big\{\d m_t\!\circ\!\d\xi\!-\!m_t\na\xi\big\}\big|_x
=\Im\,\d m_t\circ\big\{\d\xi\!-\!\na\xi\big\}\big|_x\\
&\equiv \d m_t\big(T_{t\xi(x)}E^{\nh}\big).
\end{split}\end{equation*}
The last equality on the first line follows from~\e_ref{bundleseq_e2}.
These two observations imply~\e_ref{genconn_e2b}.

\subsection{Metric-compatible connections}
\label{metricconn_subs}

\noindent
Suppose $E\!\lra\!M$ is a smooth vector bundle.
Let $g$ be a metric on $E$, i.e.
$$g\in\Ga(M;E^*\!\otimes_{\R}\!E^*) \qquad\hbox{s.t.}\qquad
g(v,w)=g(w,v),~~g(v,v)>0~~~\forall~v,w\in E_x,~v\!\neq\!0,~x\!\in\!M.$$
A connection $\na$ in $E$ is \sf{$g$-compatible} if 
$$\d\big(g(\xi,\ze)\big)=g(\na\xi,\ze)+g(\xi,\na\ze)\in \Ga(M;T^*M)
\qquad\forall~\xi,\ze\in\Ga(M;E).$$\\

\noindent
Suppose $\U$ is an open subset of $M$ and $\xi_1,\ldots,\xi_n\!\in\!\Ga(\U;E)$ is
a frame for $E$ on~$\U$.
For $i,j\!=\!1,\ldots,n$, let
$$g_{ij}=g(\xi_i,\xi_j)\in C^{\i}(\U).$$
If $\na$ is a connection in $E$ and $\th_{kl}$ is
the connection one-form for $\na$ with respect to 
the frame $\{\xi_k\}_k$, then $\na$ is $g$-compatible on~$\U$
if and only~if 
\begin{equation}\label{metriccomp_e}
\sum_{k=1}^{k=n}\big(g_{ik}\th^k_j+g_{jk}\th^k_i\big)=\d g_{ij}
\qquad\forall~i,j=1,2,\ldots,n.
\end{equation}

\subsection{Torsion-free connections}
\label{LCconn_subs}

\noindent
If $M$ is a smooth manifold, a connection $\na$ in $TM$ is \sf{torsion-free}
if
$$\na_XY-\na_YX=[X,Y].$$\\

\noindent
If $(x_1,\ldots,x_n)\!:\U\!\lra\!\R^n$ is a coordinate chart on $M$, let
$$\frac{\partial}{\partial x_1},\ldots,\frac{\partial}{\partial x_n}
\in\Ga(\U;TM)$$
be the corresponding frame for $TM$ on $\U$.
If $\na$ is a connection in~$TM$, the corresponding connection one-form~$\th$
can be written~as
$$\th_j^k=\sum_{i=1}^{i=n}\Ga_{ij}^k\d x^i, 
\qquad\hbox{where}\qquad
\na_{\partial/\partial x_i}\frac{\partial}{\partial x_j}
=\sum_{k=1}^{k=n}\Ga_{ij}^k\frac{\partial}{\partial x_k}.$$
The connection $\na$ is torsion-free on $TM|_{\U}$ if and only if
\begin{equation}\label{torsionfree_e}
\Ga_{ij}^k=\Ga_{ji}^k \qquad\forall~i,j,k=1,\ldots,n.
\end{equation}

\begin{lmm}
\label{LCconn_lmm}
If $(M,g)$ is a Riemannian manifold, there exists a unique torsion-free
$g$-compatible connection~$\na$ in~$TM$.
\end{lmm}

\noindent
{\it Proof:} (1) Suppose $\na$ and $\ti\na$ are 
torsion-free $g$-compatible connections in~$TM$.
By~\e_ref{connection_diff_e}, there exists 
\begin{gather*}
\th\in\Ga\big(M;T^*M\!\otimes_{\R}\!\Hom_{\R}(TM,TM)\big) \qquad\hbox{s.t.}\\
\ti\na_X Y-\na_X Y=\big\{\th(X)\big\}Y 
\qquad\forall~Y\!\in\!\Ga(M;TM),~X\!\in\!T_xM,~x\!\in\!M.
\end{gather*}
Since $\na$ and $\ti\na$ are torsion-free,
\begin{equation}\label{LCconn_e1}
\big\{\th(X)\big\}Y=\big\{\th(Y)\big\}X 
\qquad\forall~X,Y\in T_xM,~x\!\in\!M.
\end{equation}
Since $\na$ and $\ti\na$ are $g$-compatible,
\begin{equation}\label{LCconn_e2}
\begin{cases}
g\big(\{\th(X)\}Y,Z\big)+g\big(Y,\{\th(X)\}Z\big)=0\\
g\big(\{\th(Y)\}X,Z\big)+g\big(X,\{\th(Y)\}Z\big)=0\\
g\big(\{\th(Z)\}X,Y\big)+g\big(X,\{\th(Z)\}Y\big)=0
\end{cases}
\qquad\forall~X,Y,Z\in T_xM,~x\!\in\!M.
\end{equation}
Adding the first two equations in~\e_ref{LCconn_e2},
subtracting the third, and using~\e_ref{LCconn_e1} 
and the symmetry of~$g$, we obtain
$$2g\big(\{\th(X)\}Y,Z\big)=0 \quad\forall~X,Y,Z\in T_xM,~x\!\in\!M
\qquad\Lra\qquad \th\equiv0.$$
Thus, $\ti\na\!=\!\na$.\\

\noindent
(2) Let $(x_1,\ldots,x_n)\!:\U\!\lra\!\R^n$ be a coordinate chart on $M$.
With notation as in the paragraph preceding Lemma~\ref{LCconn_lmm},
$\na$ is $g$-compatible on $TM|_{\U}$ if and only if 
\begin{equation}\label{LCconn_e3}
\sum_{l=1}^{l=n}\big(g_{il}\Ga_{kj}^l+g_{jl}\Ga_{ki}^l\big)
=\partial_{x_k}g_{ij};
\end{equation}
see~\e_ref{metriccomp_e}.
Define a connection $\na$ in $TM|_{\U}$ by
$$\Ga_{ij}^k=\frac{1}{2}\sum_{l=1}^{l=n}g^{kl}\big(
\partial_{x_i}g_{jl}+\partial_{x_j}g_{il}-\partial_{x_l}g_{ij}\big)
\qquad\forall~i,j,k=1,\ldots,n,$$
where $g^{ij}$ is the $(i,j)$-entry of the inverse of the matrix
$(g_{ij})_{i,j=1,\ldots,n}$.
Since $g_{ij}\!=\!g_{ji}$, $\Ga_{ij}^k$ satisfies \e_ref{torsionfree_e};
a direct computation shows that $\Ga_{ij}^k$ also satisfies \e_ref{LCconn_e3}. 
Therefore, $\na$ is a torsion-free $g$-compatible connection on~$TM|_{\U}$.
In this way, we can define a  torsion-free $g$-compatible connection
on every coordinate chart.
By the uniqueness property, these connections agree on the overlaps.

\section{Complex structures}

\subsection{Complex linear connections}
\label{Cconn_subs}

\noindent
Suppose $M$ is a smooth manifold and $\pi\!:(E,\fI)\!\lra\!M$ is 
a complex vector bundle.
Similarly to Subsection~\ref{genconn_subs}, there is an exact
sequence of vector bundles 
\begin{equation}\label{Cbundleseq_e1}
0 \lra \pi^*E \lra TE \stackrel{\d\pi}{\lra} \pi^*TM\lra 0
\end{equation}
over~$E$.
If $f\!\in\!C^{\i}(M;\C)$ and $m_f\!: E\!\lra\!E$ is defined as
in~\e_ref{mf_dfn_e}, there is a commutative diagram
\begin{equation}\label{Cbundleseq_e2}
\xymatrix{0 \ar[r]& \pi^*E \ar[r] \ar[d]^{\pi^*m_f} & TE \ar[r]^{\d\pi} \ar[d]^{\d m_f} &
\pi^*TM \ar[r] \ar[d]^{\id} & 0\\
0 \ar[r]& \pi^*E \ar[r]& m_f^*TE \ar[r]^{\d\pi} & \pi^*TM \ar[r]  & 0}
\end{equation}
of bundle maps over $E$.\\

\noindent
Suppose 
$$\na\!:\Ga(M;E)\lra\Ga(M;T^*M\!\otimes_{\R}\!E)$$
is a $\C$-linear connection, i.e.
$$\na_v(\fI\xi)=\fI(\na_v\xi) \qquad\forall\,\xi\in\Ga(M;E),\,v\!\in\!TM.$$
If $\U$ is an open subset of $M$ and $\xi_1,\ldots,\xi_n\!\in\!\Ga(\U;E)$ is
a $\C$-frame for $E$ on~$\U$, then there exist 
$$\th^k_l\in \Ga(M;T^*M) \qquad\hbox{s.t.}\qquad
\na\xi_l=\sum_{k=1}^{k=n}\xi_k\th^k_l 
\equiv \sum_{k=1}^{k=n}\th^k_l\!\otimes\!\xi_k
\quad\forall~l\!=\!1,\ldots,n.$$
We will call
$$\th \equiv \big(\th^k_l\big)_{k,l=1,\ldots,n}
\in \Ga\big(\Si;T^*M\!\otimes_{\R}\!\Mat_n\C\big)$$
\sf{the complex connection one-form of $\na$ with respect to the frame $(\xi_k)_k$}.
For an arbitrary section
$$\xi=\sum_{l=1}^{l=n}f^l\xi_l\in \Ga(\U;E),$$
by~\e_ref{na_dfn_e} and $\C$-linearity of $\na$ we have
\begin{equation}\label{Cna_e3}
\na\xi=\sum_{k=1}^{k=n}\xi_k \Big(\d f^k +
\sum_{l=1}^{l=n}\th^k_lf^l\Big), \qquad\hbox{i.e.}\qquad
\na\big(\un{\xi}\cdot\un{f}^t\big)=\un\xi\cdot\big\{\d+\th\big\}\un{f}^t,
\end{equation}
where $\un\xi$ and $\un{f}$ are as \e_ref{Cna_e3b}.\\

\noindent
Let $g$ be a hermitian metric on $E$, i.e.
$$g\in\Ga\big(M;\Hom_{\C}(\bar{E}\!\otimes_{\C}\!E,\C)\big) \quad\hbox{s.t.}\quad
g(v,w)=\ov{g(w,v)},~~g(v,v)>0~~~\forall~v,w\in E_x,~v\!\neq\!0,~x\!\in\!M.$$
A $\C$-linear connection $\na$ in $E$ is \sf{$g$-compatible} if 
$$\d\big(g(\xi,\ze)\big)=g(\na\xi,\ze)+g(\xi,\na\ze)\in\Ga(M;T^*M\!\otimes_{\R}\!\C)
\qquad\forall~\xi,\ze\in\Ga(M;E).$$
With notation as in the previous paragraph, let
$$g_{ij}=g(\xi_i,\xi_j)\in C^{\i}(\U;\C) \qquad\forall~i,j\!=\!1,\ldots,n.$$
Then $\na$ is $g$-compatible on~$\U$ if and only~if 
\begin{equation}\label{Cmetriccomp_e}
\sum_{k=1}^{k=n}\big(g_{ik}\th^k_j+\bar{g}_{jk}\bar\th^k_i\big)=\d g_{ij}
\qquad\forall~i,j=1,2,\ldots,n.
\end{equation}

\subsection{Generalized $\dbar$-operators}
\label{dbar_subs}

\noindent
If $(\Si,\fj)$ is an almost complex manifold, let
$$T^*\Si^{0,1}\equiv\big\{\eta\!\in\!T^*\Si\!\otimes\!_{\R}\C\!:
\eta\circ \fj=\fI\,\eta\big\}
\qquad\hbox{and}\qquad
T^*\Si^{0,1}\equiv\big\{\eta\!\in\!T^*\Si\!\otimes\!_{\R}\C\!:
\eta\circ \fj=-\fI\,\eta\big\}$$
be the bundles of $\C$-linear and $\C$-antilinear $1$-forms on $\Si$.
If $(\Si,\fj)$ and $(M,J)$ are smooth almost complex manifolds
and $u\!:\Si\!\lra\!M$ is a smooth function, define 
\BE{dbarJjdfn_e}
\dbar_{J,\fj} u\in \Ga\big(\Si;T^*\Si^{0,1}\!\otimes_{\C}\!u^*TM\big) \qquad\hbox{by}\qquad
\dbar_{J,\fj}u=\frac{1}{2}\big(\d u+J\circ \d u\circ \fj\big).\EE
A smooth map $u\!:(\Si,\fj)\!\lra\!(M,J)$ will be called \sf{$(J,\fj)$-holomorphic} 
if $\dbar_{J,\fj}u\!=\!0$. 

\begin{dfn}
\label{dbar_dfn}
Suppose $(\Si,\fj)$ is an almost complex manifold and $\pi\!:(E,\fI)\!\lra\!\Si$ 
is a complex vector bundle.
A \sf{$\dbar$-operator on $(E,\fI)$} is a $\C$-linear map
$$\dbar\!: \Ga(\Si;E)\lra \Ga(\Si;T^*\Si^{0,1}\!\otimes_{\C}\!E)$$
such that 
\begin{equation}\label{dbar_dfn_e}
\dbar\big(f\xi)=(\dbar{f})\!\otimes\!\xi+f(\dbar\xi)
\qquad\forall~f\!\in\!C^{\i}(\Si),~\xi\!\in\!\Ga(\Si;E),
\end{equation}
where $\dbar{f}\!=\!\dbar_{\fI,\fj}{f}$ is the usual $\dbar$-operator 
on complex-valued functions.
\end{dfn}

\noindent
Similarly to Subsection~\ref{genconn_subs},
a $\dbar$-operator on $(E,\fI)$ is a first-order differential operator.
If $\U$ is an open subset of $M$ and $\xi_1,\ldots,\xi_n\!\in\!\Ga(\U;E)$ is
a $\C$-frame for $E$ on~$\U$, then there exist 
$$\th^k_l\in\Ga(\Si;T^*\Si^{0,1}) \qquad\hbox{s.t.}\qquad
\dbar\xi_l=\sum_{k=1}^{k=n}\xi_k\th^k_l \equiv \sum_{k=1}^{k=n}\th^k_l\!\otimes\!\xi_k
\quad\forall~l\!=\!1,\ldots,n.$$
We will call
$$\th \equiv \big(\th^k_l\big)_{k,l=1,\ldots,n}
\in \Ga\big(\Si;T^*\Si^{0,1}\!\otimes_{\C}\!\Mat_n\C\big)$$
\sf{the connection one-form of $\dbar$ with respect to the frame $(\xi_k)_k$}.
For an arbitrary section
$$\xi=\sum_{l=1}^{l=n}f^l\xi_l\in \Ga(\U;E),$$
by~\e_ref{dbar_dfn_e} we have
\begin{equation}\label{dbar_e3}
\dbar\xi=\sum_{k=1}^{k=n}\xi_k \Big(\dbar f^k +
\sum_{l=1}^{l=n}\th^k_lf^l\Big), \qquad\hbox{i.e.}\qquad
\dbar\big(\un{\xi}\cdot\un{f}^t\big)=\un\xi\cdot\big\{\dbar+\th\big\}\un{f}^t,
\end{equation}
where $\un\xi$ and $\un{f}$ are as in~\e_ref{Cna_e3b}.
It is immediate from~\e_ref{dbar_dfn_e} that the symbol of~$\dbar$ is given~by
$$\si_{\dbar}\!: T^*\Si\lra  \Hom_{\C}\big(E,T^*\Si^{0,1}\!\otimes_{\C}\!E\big), \qquad
\big\{\si_{\dbar}(\eta)\big\}(f)=
\big(\eta+\fI\,\eta\circ \fj\big) \otimes f.$$
In particular, $\dbar$ is an elliptic operator (i.e.~$\si_{\dbar}(\eta)$ is an isomorphism
for $\eta\!\neq\!0$) if $(\Si,\fj)$ is a Riemann surface.

\begin{lmm}\label{complexstr_lmm1}
Suppose $(\Si,\fj)$ is an almost complex manifold and $\pi\!:(E,\fI)\!\lra\!\Si$ 
is a complex vector bundle.
If 
$$\dbar\!: \Ga(\Si;E)\lra \Ga(\Si;T^*\Si^{0,1}\!\otimes_{\C}\!E)$$
is a $\dbar$-operator on $(E,\fI)$, 
there exists a unique almost complex structure $J\!=\!J_{\dbar}$
on (the total space of) $E$ such that $\pi$ is a $(\fj,J)$-holomorphic map,
the restriction of $J$ to the vertical tangent bundle $TE^{\nv}\!\approx\!\pi^*E$
agrees with~$\fI$, and 
\begin{equation}\label{complexstr_lmm1_e}
\dbar_{J,\fj}\xi=0\in\Ga(\U;T^*\Si^{0,1}\!\otimes_{\C}\!\xi^*TE)
\qquad\Llra\qquad
\dbar\xi=0\in\Ga(\U;T^*\Si^{0,1}\!\otimes_{\C}\!E)
\end{equation}
for every open subset $\U$ of $\Si$ and $\xi\!\in\!\Ga(\U;E)$.
\end{lmm}

\noindent
{\it Proof:} (1) With notation as above, define
$$\vph\!: \U\!\times\!\C^n \lra E|_{\U} \qquad\hbox{by}\qquad
\vph(x,c^1,\ldots,c^n)=\un\xi(x)\cdot\un{c}^t\equiv \sum_{k=1}^{k=n}c^k\xi_k(x)\in E_x.$$
The map $\vph$ is a trivialization of $E$ over~$\U$.
If $J\!\equiv\!J_{\dbar}$ is an almost complex structure on~$E$,
let $\ti{J}$ be the almost complex structure on $\U\!\times\!\C^n$ given~by
\begin{equation}\label{complexstr_lmm1e0}
\ti{J}_{(x,\un{c})}=\big\{\d_{(x,\un{c})}\vph\big\}^{-1}
\circ J_{\vph(x,\un{c})}\circ \d_{(x,\un{c})}\vph
\qquad\forall~(x,\un{c})\in \U\!\times\!\C^n.
\end{equation}
The almost complex structure $J$ restricts to $\fI$ on $TE^{\nv}$ if and only if
\begin{equation}\label{complexstr_lmm1e1}
\ti{J}_{(x,\un{c})}w=\fI w \in 
T_{\un{c}}\C^n\subset T_{(x,\un{c})}(\U\!\times\!\C^n) 
\qquad\forall~w\in T_{\un{c}}\C^n.
\end{equation}
If $J$ restricts to $\fI$ on $TE^{\nv}$, the projection $\pi$ is $(\fj,J)$-holomorphic
on $E|_{\U}$ if and only if there exists 
\begin{gather}
\ti{J}^{\nv\nh}\in \Ga\big(\U\!\times\!\C^n;\Hom_{\R}(\pi_{\U}^*T\U,\pi_{\C^n}^*T\C^n)\big)
\qquad\hbox{s.t.}\notag\\
\label{complexstr_lmm1e2}
\ti{J}_{(x,\un{c})}w=\fj_x w+\ti{J}_{(x,\un{c})}^{\nv\nh}w
\qquad\forall~w\in T_x\U\subset T_{(x,\un{c})}(\U\!\times\!\C^n).
\end{gather}
If $\xi\!\in\!\Ga(\U;E)$, let 
$$\ti\xi\equiv \vph^{-1}\circ \xi\equiv \big(\id_{\U},\un{f}\big),
\qquad\hbox{where}\qquad \un{f}\in C^{\i}(\U;\C^n).$$
By \e_ref{complexstr_lmm1e0}-\e_ref{complexstr_lmm1e2},
\begin{equation}\label{complexstr_lmm1e3}\begin{split}
2\,\dbar_{J,\fj}\xi\big|_x = \d_{\ti\xi(x)}\vph\circ 2\dbar_{\ti{J},\fj}\ti\xi\big|_x
&=\d_{\ti\xi(x)}\vph\circ\big\{\big(\Id_{T_x\U},\d_x\un{f}\big)
+\ti{J}_{\ti\xi(x)}\circ\big(\Id_{T_x\U},\d_x\un{f}\big)\circ\fj_x\big\}\\
&=\d_{\ti\xi(x)}\vph\circ\big(0,2\,\dbar{f}|_x+\ti{J}^{\nv\nh}_{\ti\xi(x)}\circ \fj_x\big).
\end{split}\end{equation}
On the other hand, by~\e_ref{dbar_e3},
\begin{equation}\label{complexstr_lmm1e4}\begin{split}
\dbar\xi|_x=\dbar(\un\xi\cdot{f}^t\big)\big|_x
&=\un\xi(x)\cdot\big\{\dbar\!+\!\th\}{f}^t\big|_x\\
&=\vph\big(\dbar{f}|_x+\th_x\cdot f(x)^t\big).
\end{split}\end{equation}
By~\e_ref{complexstr_lmm1e3} and \e_ref{complexstr_lmm1e4},
the property~\e_ref{complexstr_lmm1_e} is satisfied for all $\xi\!\in\!\Ga(\U;E)$ if
and only if
$$\ti{J}^{\nv\nh}_{(x,\un{c})}
=2\big(\th_x\cdot\un{c}^t\big)\circ(-\fj_x) =2\fI\,\th_x\cdot\un{c}^t
\qquad\forall~(x,\un{c})\in\U\!\times\!\C^n.$$
In summary, the almost complex structure $J\!=\!J_{\dbar}$ on $E$ has
the three desired properties if and only if for every trivialization of $E$ over 
an open subset~$\U$ of~$\Si$
\begin{gather}\label{complexstr_lmm1e6}
\ti{J}_{(x,\un{c})}\big(w_1,w_2\big)
=\big(\fj_x w_1,\fI w_2+2\fI\th_x(w_1)\cdot\un{c}^t\big)\\
\qquad\forall~(x,\un{c})\in\U\!\times\!\C^n,~
(w_1,w_2)\in T_x\U\!\oplus\!T_{\un{c}}\C^n=T_{(x,\un{c})}(\U\!\times\!\C^n), \notag
\end{gather}
where $\ti{J}$ is the almost complex structure on $\U\!\times\!\C^n$ induced
by $J$ via the trivialization and $\th$ is the connection-one form corresponding	
to $\dbar$ with respect to the frame inducing the trivialization.\\

\noindent
(2) By~\e_ref{complexstr_lmm1e6}, there exists at most one almost complex structure 
$J$ satisfying the three properties.
Conversely, \e_ref{complexstr_lmm1e6} determines such an almost complex structure on~$E$.
Since
\begin{equation*}\begin{split}
\ti{J}_{(x,\un{c})}^2\big(w_1,w_2\big)
=\ti{J}_{(x,\un{c})}\big(\fj w_1,\fI w_2+2\fI\th_x(w_1)\cdot\un{c}^t\big)
&=\big(\fj^2 w_1,\fI\big(\fI w_2+2\fI\th_x(w_1)\cdot\un{c}^t\big)+
2\fI\th_x(\fj w_1)\cdot\un{c}^t\big)\\
&=-(w_1,w_2),
\end{split}\end{equation*}
$\ti{J}$ is indeed an almost complex structure on~$E$.
The almost complex structure induced by $\ti{J}$ on~$E|_{\U}$ satisfies
the three properties by part~(a).
By the uniqueness property, the almost complex structures on $E$ induced by
the different trivializations agree on the overlaps.
Therefore, they define an almost complex structure $J\!=\!J_{\dbar}$ 
on the total space of $E$ with the desired properties.

\subsection{Connections and $\dbar$-operators}

\noindent
Suppose $(\Si,\fj)$ is an almost complex manifold, $\pi\!:(E,\fI)\!\lra\!\Si$ 
is a complex vector bundle, and
$$\dbar\!: \Ga(\Si;E)\lra \Ga(\Si;T^*\Si^{0,1}\!\otimes_{\C}\!E)$$
is a $\dbar$-operator on $(E,\fI)$.
A $\C$-linear connection $\na$ in $(E,\fI)$ is \sf{$\dbar$-compatible}
if
\begin{equation}\label{dbarna_e}
\dbar\xi=\dbar_{\na}\xi \equiv
\frac{1}{2}\big(\na\xi+\fI\na\xi\circ\fj\big)
\qquad\forall~\xi\!\in\!\Ga(M;\Si).
\end{equation}

\begin{lmm}
\label{dbarconn_lmm}
Suppose $(\Si,\fj)$ is an almost complex manifold, $\pi\!:(E,\fI)\!\lra\!\Si$ 
is a complex vector bundle,
$$\dbar\!: \Ga(\Si;E)\lra \Ga(\Si;T^*\Si^{0,1}\!\otimes_{\C}\!E)$$
is a $\dbar$-operator on $(E,\fI)$, and $J_{\dbar}$ is
the complex structure in the vector bundle $TE\!\lra\!E$
provided by Lemma~\ref{complexstr_lmm1}.
A $\C$-linear connection $\na$ in $(E,\fI)$ is $\dbar$-compatible
if and only if the splitting~\e_ref{genconn_e1} determined by~$\na$
respects the complex structures.
\end{lmm}

\noindent
{\it Proof:} Since $J_{\dbar}\!=\!\pi^*\fI$ on $\pi^*E\!\subset\!TE$, 
the splitting~\e_ref{genconn_e1} determined by~$\na$ respects the complex structures
if and only~if
$$J_{\dbar}|_v\circ\big\{\d\xi-\na\xi\big\}\big|_x
=\big\{\d\xi-\na\xi\big\}\big|_x\circ\fj_x\!: T_x\Si\lra T_vE$$
for all $x\!\in\!\Si$, $v\!\in\!E_x$, and $\xi\!\in\!\Ga(\Si;E)$ 
such that $\xi(x)\!=\!0$; see the proof of Lemma~\ref{genconn_lmm}.
This identity is equivalent to
\begin{equation}\label{dbarnconn_e1}
\dbar_{J_{\dbar},\fj}\xi = \dbar_{\na}\xi \qquad\forall~\xi\in\Ga(\Si;E).
\end{equation}
On the other hand, by the proof of Lemma~\ref{complexstr_lmm1}, 
\begin{equation}\label{dbarnconn_e2}
\dbar_{J_{\dbar},\fj}\xi = \dbar\xi \qquad\forall~\xi\in\Ga(\Si;E);
\end{equation}
see~\e_ref{complexstr_lmm1e3}-\e_ref{complexstr_lmm1e6}.
The lemma follows immediately from~\e_ref{dbarnconn_e1} and~\e_ref{dbarnconn_e2}.

\subsection{Holomorphic vector bundles}

\noindent
Let $(\Si,\fj)$ be a complex manifold.
A \sf{holomorphic vector bundle} $(E,\fI)$ on $(\Si,\fj)$
is a complex vector bundle with a collection of trivializations
that overlap holomorphically.\\

\noindent
A collection of holomorphically overlapping trivializations of $(E,i)$
determines a holomorphic structure~$J$ on the total space of~$E$ and 
a $\dbar$-operator
$$\dbar\!: \Ga(\Si;E)\lra \Ga(\Si;T^*\Si^{0,1}\!\otimes_{\C}\!E).$$
The latter is defined as follows.
If $\xi_1,\ldots,\xi_n$ is a holomorphic complex frame for $E$
over an open subset $\U$ of~$M$, then
$$\dbar\sum_{k=1}^{k=n}f^k\xi_k 
=\sum_{k=1}^{k=n}\dbar f^k\!\otimes\!\xi_k
\qquad\forall~f^1,\ldots,f^k\in C^{\i}(\U;\C).$$
In particular, for all $\xi\!\in\!\Ga(M;E)$
$$\dbar_{J,\fj}\xi=0 \qquad\Llra\qquad \dbar\xi=0.$$
Thus, $J\!=\!J_{\dbar}$; see Lemma~\ref{complexstr_lmm1}.

\begin{lmm}
\label{complexstr_lmm2}
Suppose $(\Si,\fj)$ is a Riemann surface and $\pi\!:(E,\fI)\!\lra\!\Si$ 
is a complex vector bundle.
If 
$$\dbar\!: \Ga(\Si;E)\lra \Ga(\Si;T^*\Si^{0,1}\!\otimes_{\C}\!E)$$
is a $\dbar$-operator on $(E,\fI)$,
the almost complex structure $J\!=\!J_{\dbar}$ on $E$ is integrable.
With this complex structure, $\pi\!:E\!\lra\!\Si$ is a holomorphic vector bundle
and $\dbar$ is the corresponding $\dbar$-operator.
\end{lmm}

\noindent
{\it Proof:} By~\e_ref{complexstr_lmm1_e}, it is sufficient to show that there exists a 
$(J,\fj)$-holomorphic local section through every point $v\!\in\!E$,
i.e.~there exist a neighborhood $\U$ of $x\!\equiv\!\pi(v)$ in $\Si$
and $\xi\!\in\!\Ga(\U;E)$ such that 
$$\xi(x)=v \qquad\hbox{and}\qquad \dbar_{J,\fj}\xi=0.$$
By Lemma~\ref{complexstr_lmm1} and~\e_ref{complexstr_lmm1e4},
this is equivalent to showing that the equation
\begin{equation}\label{complexstr_lmm2e1}
\Big\{\dbar+\th\Big\}f^t=0, \qquad f(x)=v, \qquad f\in C^{\i}(\U;\C^n),
\end{equation}
has a solution for every $v\!\in\!\C^n$. 
We can assume that $\U$ is a small disk contained in $S^2$.
Let 
$$\eta\!:S^2\lra[0,1]$$ 
be a smooth function supported in $\U$ and such that $\eta\!\equiv\!1$
on a neighborhood of~$x$.
Then, 
$$\eta\th\in  \Ga(S^2;(T^*S^2)^{0,1}\!\otimes_{\C}\!\Mat_n\C).$$
Choose $p\!>\!2$.
The operator
$$\Th: L^p_1(S^2;\C^n) \lra 
L^p\big(S^2;(T^*S^2)^{0,1}\!\otimes_{\C}\!\C^n\big)\oplus\C^n,
\qquad
\Th(f)=\big(\dbar_{\fI,\fj}f,f(x)\big),$$
is surjective.
If $\eta$ has sufficiently small support, so is the operator
$$\Th_{\eta}: L^p_1(S^2;\C^n) \lra 
L^p\big(S^2;(T^*S^2)^{0,1}\!\otimes_{\C}\!\C^n\big)\oplus\C^n,
\qquad
\Th_{\eta}(f)=\big(\{\dbar_{\fI,\fj}\!+\!\eta\th\}f,f(x)\big).$$
Then, the restriction of $\Th_{\eta}^{-1}(0,v)$ to a neighborhood of $x$
on which $\eta\!\equiv\!1$ is a solution of~\e_ref{complexstr_lmm2e1}.
By elliptic regularity, $\Th_{\eta}^{-1}(0,v)\!\in\!C^{\i}(S^2;\C^n)$.

\subsection{Deformations of almost complex submanifolds}
\label{complexmfld_subs}

\noindent
If $(M,J)$ is a complex manifold, holomorphic coordinate charts on $(M,J)$ 
determine a holomorphic structure in the vector bundle $(TM,\fI)\!\lra\!M$.
If $(\Si,\fj)\!\subset\!(M,J)$ is a complex submanifold,
holomorphic coordinate charts on $\Si$ can be extended to holomorphic 
coordinate charts on~$M$. 
Thus, the holomorphic structure in $T\Si\!\lra\!\Si$ induced 
from~$(\Si,\fj)$ is the restriction of the holomorphic structure in~$TM|_{\Si}$.
It follows that 
$$\dbar_M=\dbar_{\Si}\!: \Ga(\Si;T\Si)\lra 
\Ga(\Si;T^*\Si^{0,1}\!\otimes_{\C}\!T\Si)\subset 
\Ga\big(\Si;T^*\Si^{0,1}\!\otimes_{\C}\!TM|_{\Si}\big),$$
where $\dbar_M$ and $\dbar_{\Si}$ are the $\dbar$-operators in $TM|_{\Si}$ 
and $T\Si$ induced from the holomorphic structures in~$\Si$ and~$M$.
Therefore, $\dbar_M$ descends to a $\dbar$-operator  on the quotient
$$\dbar\!:\Ga(\Si;\N_M\Si)
=\Ga(\Si;TM|_{\Si})\big/\Ga(\Si;T\Si)
\lra\Ga(\Si;T^*\Si^{0,1}\!\otimes_{\C}\!\N_M\Si),$$
where
$$\N_M\Si\equiv TM|_{\Si}\big/T\Si\lra\Si$$
is the normal bundle of $\Si$ in $M$.
This vector bundle inherits a holomorphic structure from that of $TM|_{\Si}$
and~$\Si$.
The above $\dbar$-operator on $\N_M$ is the $\dbar$-operator corresponding
to this induced holomorphic structure on~$\N_M\Si$.\\

\noindent
Suppose $(M,J)$ is an almost complex manifold and
$(\Si,\fj)\!\subset\!(M,J)$ is an almost complex submanifold.
Let $\na$ be a torsion-free connection in~$TM$.
Define 
\begin{gather}
D_{J;\Si}\!:\Ga(\Si;TM|_{\Si})\lra \Ga(\Si;T^*\Si^{0,1}\!\otimes_{\C}\!TM|_{\Si})
\qquad\hbox{by}\notag\\
\label{DJSi_e}
D_{J;\Si}\xi=\frac{1}{2}\big(\na\xi+J\circ\na\xi\circ\fj\big)
-\frac{1}{2}J\circ\na_{\xi}J\!: T\Si\lra TM|_{\Si}.
\end{gather}
If $\na$ is the Levi-Civita connection (the connection of Lemma~\ref{LCconn_lmm})
for a $J$-compatible metric on $M$
(and $\Si$ is a Riemann surface), then $D_{J;\Si}$ is the linearization
of the $\dbar_J$-operator at the inclusion map $\io\!:\Si\!\lra\!M$;
see \cite[Proposition~3.1.1]{MS}.\\
 
\noindent
In fact, $D_{J;\Si}$ is independent of the choice of a torsion-free connection in $TM$.
Let
\begin{equation}\label{CRoper_e1}
\ti\na=\na+\th, \qquad\th\in\Ga\big(M;T^*M\!\otimes_{\R}\!\Hom_{\R}(TM,TM)\big),
\end{equation}
be another torsion-free connection; see~\e_ref{connection_diff_e}.
Since $\ti\na$ and $\na$ are torsion-free connections,
\begin{equation}\label{CRoper_e3}
\big\{\th(X)\big\}Y=\big\{\th(Y)\big\}X \qquad\forall X,Y\!\in\!T_xM,~x\!\in\!M. 
\end{equation}
If $x\!\in\!M$ and $X,Y\!\in\!\Ga(M;TM)$,
\begin{gather}
\big\{\na_YJ\big\}X=\na_Y(JX) - J\na_YX\,, \quad
\big\{\ti\na_YJ\big\}X=\ti\na_Y(JX) - J\ti\na_YX\qquad\Lra\notag\\
\label{CRoper_e5a}
\big\{\ti\na_YJ\big\}X-\big\{\na_YJ\big\}X
=\big\{\th(Y)\big\}(JX)-J\big\{\th(Y)\big\}X
=\big\{\th(JX)\big\}Y-J\big\{\th(X)\big\}Y
\end{gather}
by~\e_ref{CRoper_e1} and~\e_ref{CRoper_e3}.
On the other hand, by~\e_ref{CRoper_e1} for all $X\!\in\!T\Si$
and $\xi\!\in\!\Ga(\Si;TM|_{\Si})$,
\begin{equation}\label{CRoper_e5b}\begin{split}
\big\{\ti\na\xi+J\circ\ti\na\xi\circ\fj\big\}(X)
-\big\{\na\xi+J\circ\na\xi\circ\fj\big\}(X)
=\big\{\th(X)\big\}\xi+J\big\{\th(\fj X)\big\}\xi\\
=J\big(\big\{\th(JX)\big\}\xi-J\big\{\th(X)\big\}\xi\big),
\end{split}\end{equation}
since $\fj\!=\!J|_{T\Si}$ and $J^2\!=\!-\Id$.
By~\e_ref{CRoper_e5a} and~\e_ref{CRoper_e5b}, 
$D_{J,\Si}$ is independent of the choice of torsion-free connection~$\na$.\\

\noindent
Since any torsion-free connection on $\Si$ extends to
a torsion-free connection on $M$, the above observation implies that 
\BE{DJSIrestr_e} D_{J;\Si}\!: \Ga(\Si;T\Si)\lra \Ga(\Si;T^*\Si^{0,1}\!\otimes_{\C}\!T\Si)
\subset \Ga(\Si;T^*\Si^{0,1}\!\otimes_{\C}\!TM|_{\Si}).\EE
Thus, an almost complex submanifold $(\Si,\fj)$ of an almost complex manifold $(M,J)$
induces a well-defined generalized Cauchy-Riemann operator\footnote{see Section~\ref{ellipticest_subs}} 
on the normal bundle of $\Si$ in~$M$,
$$D_{J;\Si}^{\N}\!: \Ga(\Si;\N_M\Si)\lra \Ga(\Si;T^*\Si^{0,1}\!\otimes_{\C}\!\N_M\Si),
\qquad D_{J;\Si}^{\N}\big(\pi(\xi)\big)= \pi\big(D_{J;\Si}(\xi)\big)~~
\forall\,\xi\!\in\!\Ga(\Si;TM|_{\Si}),$$
where $\pi\!:TM|_{\Si}\!\lra\!\N_M\Si$ is the quotient projection map.
The $\C$-linear part of~$D_{J;\Si}^{\N}$ determines a $\dbar$-operator
on the normal bundle of $\Si$ in~$M$:
\begin{gather*}\label{dbarN_e}
\dbar_{J;\Si}^{\N}\!: \Ga(\Si;\N_M\Si)\lra \Ga(\Si;T^*\Si^{0,1}\!\otimes_{\C}\!\N_M\Si),\\
\dbar_{J;\Si}^{\N}(\xi)
=\frac{1}{2}\big(D_{J;\Si}^{\N}(\xi) -J D_{J;\Si}^{\N}(J\xi)\big)
~~\forall\,\xi\!\in\!\Ga(\Si;\N_M\Si).\notag
\end{gather*}
Both operators are determined by the almost complex submanifold $(\Si,\fj)$ 
of the almost complex manifold $(M,J)$ only and are independent of the choice of
torsion-free connection~$\na$ in~\e_ref{DJSi_e}.\\

\noindent
Any connection $\na$ in $TM$ induces a $J$-linear connection in $TM$ by
\BE{Jconn_e}\na^J_X\xi=\na_X\xi-\frac{1}{2}J(\na_XJ)\xi
\qquad\forall\,X\!\in\!TM,\,\xi\!\in\!\Ga(M;TM).\EE
If $\na$ is as in \e_ref{DJSi_e}, 
\BE{DJSi_e2} 
\big\{D_{J;\Si}\xi\big\}(X)=\big\{\dbar_{\na^J}\xi\big\}(X)+A_J(X,\xi)
-\frac{1}{4}\big\{(\na_{J\xi}J)+J(\na_{\xi}J)\big\}(X)\EE
for all $\xi\!\in\!\Ga(\Si;TM|_{\Si})$ and $X\!\in\!T\Si$,
where $A_J$ is the Nijenhuis tensor of $J$:
\BE{Nijendfn_e} A_J(\xi_1,\xi_2)=\frac{1}{4}\Big([\xi_1,\xi_2]+J[\xi_1,J\xi_2]+J[J\xi_1,\xi_2]
-[J\xi_1,J\xi_2]\Big)
\qquad\forall~\xi_1,\xi_2\in\Ga(M;TM).\EE
Since the sum of the terms in the curly brackets in \e_ref{DJSi_e2} is $\C$-linear in $\xi$,
while the Nijenhuis tensor  is $\C$-antilinear,
the $\C$-linear operator 
\BE{dbarrest_e}
\Ga(\Si;TM|_{\Si})\lra \Ga(\Si;T^*\Si^{0,1}\!\otimes_{\C}\!TM|_{\Si}), \qquad
\xi\lra \dbar_{\na^J}(\xi)-\frac{1}{4}\big\{(\na_{J\xi}J)+J(\na_{\xi}J)\big\},\EE
takes $\Ga(\Si;T\Si)$ to $\Ga(\Si;T^*\Si^{0,1}\!\otimes_{\C}\!T\Si)$
by \e_ref{DJSIrestr_e}.
Thus, it induces a $\dbar$-operator  on $\N_M\Si$ 
and this induced operator is~$\dbar_{J;\Si}^{\N}$.
If the image of the homomorphism
$$TM\lra T^*\Si^{0,1}\otimes_{\C}TM|_{\Si}\,, \qquad
\xi\lra \na_{\xi}J-J\na_{J\xi}J\,,$$
is contained in $T^*\Si^{0,1}\!\otimes_{\C}\!T\Si$,
then $\dbar_{\na^J}$ preserves~$T\Si$ and induces a $\dbar$-operator
$\dbar_{\na^J}^{\N}$ on $\N_M\Si$ with $\dbar_{\na^J}^{\N}\!=\!\dbar_{J;\Si}^{\N}$.
In this case,
$$D_{J;\Si}^{\N}\big(\pi(\xi)\big)=\pi\big(\dbar_{\na^J}\xi+A_J(\cdot,\xi)\big)\!: 
T\Si\lra \N_M{\Si} \qquad\forall~\xi\!\in\!\Ga(\Si;TM|_{\Si}).$$
This is the case in particular if $J$ is compatible with a symplectic form~$\om$ on $M$ and $\na$ 
is the Levi-Civita connection for the metric $g(\cdot,\cdot)\!=\!\om(\cdot,J\cdot)$,
as the sum in the curly brackets in~\e_ref{DJSi_e2} then vanishes by \cite[(C.7.5)]{MS}.\\

\noindent
It is immediate that $A_J$ takes $T\Si\!\otimes_{\R}\!T\Si$ to $T\Si$ and thus induces
a bundle homomorphism 
$$A_J^{\N}\!: T\Si\otimes_{\R}\N_M\Si\lra\N_M\Si\,.$$
If $\ze$ is any vector field on $M$ such that $\ze(x)\!=\!X\!\in\!T_x\Si$ for some $x\!\in\!\Si$,
then
\BE{DJSIintr_e}\begin{split}
\big\{D_{J;\Si}\xi\}(X)&=\frac{1}{2}\big([\ze,\xi]+J[J\ze,\xi]\big)\big|_x,\\
\Big\{\dbar_{\na^J}(\xi)-\frac{1}{4}\big((\na_{J\xi}J)+J(\na_{\xi}J)\big)\Big\}(X)
&=\frac{1}{4}\big([\ze,\xi]+J[J\ze,\xi]-J[\ze,J\xi]+[J\ze,J\xi]\big)\big|_x,
\end{split}\EE
since $\na$ is torsion-free.\footnote{Since LHS and RHS of these identities depend only
$\xi$ and $X\!=\!\ze(x)$, and not on~$\ze$, it is sufficient to verify them under the assumption
that $\na\ze|_x\!=\!0$.}
These two identities immediately imply that the operators~\e_ref{DJSi_e} and~\e_ref{dbarrest_e}
preserve $T\Si\!\subset\!TM|_{\Si}$ and thus induce operators 
$$\Ga(\Si;\N_M\Si)\lra \Ga(\Si;T^*\Si^{0,1}\!\otimes_{\C}\!\N_M\Si)$$
as claimed above.\\ 

\noindent
If $g$ is a $J$-compatible metric on $TM|_{\Si}$ and 
$\pi^{\perp}\!:TM|_{\Si}\!\lra\!T\Si^{\perp}$ is the projection to 
the $g$-orthogonal complement of $T\Si$ in~$TM|_{\Si}$, the composition~$\na^{\perp}$
$$\Ga(\Si;T\Si^{\perp}) \hookrightarrow\Ga\big(\Si;TM|_{\Si}\big)
\stackrel{\na^J}{\lra} \Ga\big(\Si;T^*\Si\!\otimes\!_{\R}TM|_{\Si}\big)
\stackrel{\pi^{\perp}}{\lra} \Ga\big(\Si;T^*\Si\!\otimes\!_{\R}T\Si^{\perp}\big),$$
with $\na^J$ as in~\e_ref{Jconn_e}, is a $g$-compatible $J$-linear connection
in $T\Si^{\perp}$.
Via the isomorphism\linebreak $\pi\!: T\Si^{\perp}\!\lra\!\N_M\Si$, 
it induces a $J$-linear connection~$\na^{\N}$ in $\N_M\Si$ which is compatible 
with the metric~$g^{\N}$ induced via this isomorphism from $g|_{T\Si^{\perp}}$.
If the image of the homomorphism
\BE{nahomom_e}T\Si^{\perp}\lra T^*\Si^{0,1}\otimes_{\C}TM|_{\Si}\,, \qquad
\xi\lra \na_{\xi}J-J\na_{J\xi}J\,,\EE
is contained in $T^*\Si^{0,1}\!\otimes_{\C}T\Si$, then $\dbar_{\na^{\N}}\!=\!\dbar_{J;\Si}^{\N}$
and so
$$D_{J;\Si}^{\N}\big(\pi(\xi)\big)=\pi\big(\dbar_{\na^{\perp}}\xi+A_J(\cdot,\xi)\big)\!: 
T\Si\lra \N_M{\Si}
\qquad\forall~\xi\!\in\!\Ga(\Si;T\Si^{\perp}).$$
This is the case if $\Si$ is a divisor in $M$, i.e.~$\rk_{\C}\N=1$,
since $(\na_{\ze}J)\xi$ is $g$-orthogonal to $\xi$ and $J\xi$ for all $\xi,\ze\!\in\!T_xM$
and $x\!\in\!M$ by \cite[(C.7.1)]{MS}.
This is also the case if $J$ is compatible with a symplectic form~$\om$ on $M$ and $g(\cdot,\cdot)\!=\!\om(\cdot,J\cdot)$,
as the homomorphism~\e_ref{nahomom_e} is then trivial by~\cite[(C.7.5)]{MS}.

\section{Riemannian geometry estimates}
\label{riem_sec}

\noindent
This section is based on \cite[Chapter~1]{C} and \cite[Section~3]{F}
and culminates in a Poincare lemma for closed curves in Proposition~\ref{poincare_prp}
and an expansion for the $\dbar$-operator in Proposition~\ref{dbar_prp}.
If $u\!:\Si\!\lra\!M$ is a smooth map between smooth manifolds
and $E\!\lra\!M$ is a smooth vector bundle, let
$$\Ga(u;E)=\Ga(\Si;u^*E), \qquad
\Ga^1(u;E)=\Ga(\Si;T^*\Si\!\otimes_{\R}\!u^*E).$$
We denote the subspace of compactly supported sections
in $\Ga(u;E)$ by~$\Ga_c(u;E)$.\\

\noindent
An \sf{exponential-like map on a smooth manifold} $M$ is a smooth map
$\exp\!:TM\!\lra\!M$ such that $\exp|_M\!=\!\id_M$  and
$$\d_x\exp=\big(\id_{T_xM}~\id_{T_xM}\big)\!:
T_x(TM)=T_xM\oplus T_xM\lra T_xM \qquad\forall~x\!\in\!M,$$
where the second equality is the canonical splitting of $T_x(TM)$
into the horizontal and vertical tangent space along the zero section.
Any connection $\na$ in $TM$ gives rise to a smooth map
$\exp^{\na}\!:W\!\lra\!M$ from some neighborhood $W$ of the zero section~$M$ in~$TM$;
see \cite[Section~1.3]{C}.
If $\eta\!:TM\!\lra\!\R$ is a smooth function which equals~1 on a neighborhood of
$M$ in $TM$ and~0 outside of~$W$, then
$$\exp\!: TM\lra M, \qquad v\lra \exp^{\na}\big(\eta(v)v\big),$$ 
is an exponential-like map.
If $M$ is compact, then $W$ can be taken to be all of~$TM$ and $\exp\!=\!\exp^{\na}$.\\

\noindent
If $(M,g,\exp)$ is a Riemannian manifold with an exponential-like map and $x\!\in\!M$, 
let $r_{\exp}(x)\!\in\!\R^+$ be the supremum of the numbers $r\!\in\!\R$ such that 
the restriction 
$$\exp\!: \big\{v\!\in\!T_xM\!:\,|v|\!<\!r\big\}\lra M$$
is a diffeomorphism onto an open subset of $M$.
Set
$$r_{\exp}^g(x)=\inf\big\{d_g(x,\exp(v))\!:\,v\!\in\!T_xM,\,|v|\!=\!r_{\exp}(x)\big\}\in\R^+,$$
where $d_g$ is the metric on $M$ induced by~$g$.
If $K\!\subset\!M$, let
$$r_{\exp}^g(K)=\inf_{x\in K}r_{\exp}^g(x);$$
this number is positive if $\bar{K}\!\subset\!M$ is compact.

\subsection{Parallel transport}
\label{pt_subs}

\noindent
Let $(E,\lr{,},\na)\!\lra\!M$ be a vector bundle, real or complex, 
with an inner-product $\lr{,}$ and a metric-compatible connection~$\na$.
If $\al\!:(a,b)\!\lra\!M$ is a piecewise smooth curve, denote by 
$$\Pi_{\al}\!: E_{\al(a)} \lra E_{\al(b)}$$ 
the parallel-transport map along $\al$ with respect to  the connection~$\na$.
If $\exp\!:TM\!\lra\!M$ is an exponential-like map,
$x\!\in\!M$, and $v\!\in\!T_xM$, let 
$$\Pi_v\!:E_x\lra E_{\exp(v)}$$
be the parallel transport along the curve 
$$\ga_v\!:[0,1]\lra M, \qquad \ga_v(t)=\exp(tv).$$
If $u\!:[a,b]\!\times\![c,d]\!\lra\!M$ is a smooth map, let 
$$\Pi_{\partial u}: E_{u(a,c)}\lra E_{u(a,c)}$$ 
be  the parallel transport along $u$ restricted to the boundary
of the rectangle traversed in the positive direction.
If $u\!:\Si\!\lra\!M$ is any smooth map,  $\na$ induces 
a connection
$$\na^u\!: \Ga(u;E)\lra  \Ga^1(u;E)$$
in the vector bundle $u^*E\!\lra\!\Si$. 
If $\al$ is a smooth curve as above and $\ze\!\in\!\Ga(\al;E)$,
let
$$\frac{D}{\d t}\ze=\na^{\al}_{\partial_t}\ze\in \Ga(\al;E),$$
where $\partial_t$ is the standard unit vector field on $\R$.

\begin{lmm}\label{ptrec_lmm}
If $(M,g)$ is a Riemannian manifold
and $(E,\lr{,},\na)$ is a normed vector bundle with connection over~$M$,
for every compact subset $K\!\subset\!M$ there exists $C_K\!\in\!\R^+$ 
such that for every smooth map $u\!:[a,b]\!\times\![c,d]\!\lra\!M$ 
with $\Im\,u\!\subset\!K$
$$|\Pi_{\partial u}-\I|\le C_K\int_c^d\!\!\int_a^b\!|u_s||u_t|\d s\d t,$$
where the norm of 
$(\Pi_{\partial u}\!-\!\I)\!\in\!\End(E_{u(a,c)})$ 
is computed with respect to the inner-product in $E_{u(a,c)}$.
\end{lmm}

\begin{figure}
\begin{pspicture}(0,-2)(10,1.5)
\psset{unit=.4cm}
\psline[linewidth=.02]{->}(14.5,-4)(26,-4)
\psline[linewidth=.02]{->}(15,-4.5)(15,2)
\psline[linewidth=.04](17,-3)(25,-3)\psline[linewidth=.04](17,1)(25,1)
\psline[linewidth=.04](17,-3)(17,1)\psline[linewidth=.04](25,-3)(25,1)
\pscircle*(17,-4){.1}\rput(17,-4.5){\sm{$a$}}
\pscircle*(25,-4){.1}\rput(25,-4.5){\sm{$b$}}
\pscircle*(15,-3){.1}\rput(14.5,-3){\sm{$c$}}
\pscircle*(15,1){.1}\rput(14.5,1){\sm{$d$}}
\pscircle*(17,-3){.1}\rput(16.3,-3){\sm{$v_i$}}
\psline[linewidth=.06]{->}(17,-3)(17,-.5)\rput(16.3,-1.2){\sm{$\xi_i$}}
\psline[linewidth=.02]{->}(17,-2.2)(22,-2.2)
\psline[linewidth=.02]{->}(17,-1.4)(22,-1.4)\rput(22.8,-1.8){\sm{$\ze_i$}}
\end{pspicture}
\caption{Extending a basis $\{v_i\}$ for $E_{u(a,c)}$ to 
a frame $\{\ze_i\}$ over $[a,b]\!\times\![c,d]$} 
\label{rectangle_fig}
\end{figure}
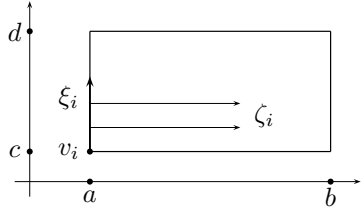

\noindent
{\it Proof:} (1) Choose an orthonormal frame $\{v_i\}$ for $E_{u(a,c)}$.
Extend each $v_i$ to
$$\xi_i\in \Ga\big(u|_{a\times[c,d]};E\big)$$
by parallel-transporting along the curve $t\!\lra\!u(a,t)$
and then to $\ze_i\!\in\!\Ga(u;E)$ by parallel-transporting $\xi_i(a,t)$ 
along the curve $s\!\lra\!u(s,t)$; see Figure~\ref{rectangle_fig}.
By construction, 
$$\frac{D}{\d s}\ze_i=0\in\Ga(u;E).$$ 
Let $A$ be the matrix-valued function on $[a,b]\!\times\![c,d]$ such~that
\BE{ptrec_e0}
\frac{D}{\d t}\ze_i\Big|_{(s,t)}=\sum_{l=1}^{l=k}A_{il}(s,t)\ze_l(s,t),\EE
where $k$ is the rank of $E$.
Note that $A_{ij}(a,t)=0$ and
\BE{ptrec_e1}\begin{split}
\blr{\cR_{\na}(u_s,u_t)\ze_i,\ze_j}
&=\bblr{\frac{D}{\d s}\frac{D}{\d t}\ze_i-
 \frac{D}{\d t}\frac{D}{\d s}\ze_i,\ze_j}
=\sum_{l=1}^{l=k}\bblr{\bigg(\frac{\partial}{\partial s}A_{il}\bigg)\ze_l,\ze_j}
=\frac{\partial}{\partial s}A_{ij}\,,
\end{split}\EE
where $\cR_{\na}$ is the curvature tensor of the connection of~$\na$. 
Since $K$ is compact and the image of~$u$ is contained in~$K$, it follows that 
\BE{ptrec_e2} |A_{ij}(b,t)|\le C_K
\int_a^b\!|u_s|_{(s,t)}|u_t|_{(s,t)}\d s. \EE

\noindent
(2) The parallel transport of $\ze_i$ along the curves
$$\tau\lra u(\tau,c),\quad\tau\lra u(\tau,d),\quad\tau\lra u(a,\tau)$$
is $\ze_i$ itself. Thus, it remains to estimate the parallel
transport of each $\ze_i$ along the curve \hbox{$\tau\!\lra\!u(b,\tau)$.}
Let $h_{ij}$ be the $\SO_k$-valued function ($\U_k$-valued function if $E$ is complex)
on $[c,d]$ such that
$$h(c)=\I,\qquad
\sum_{j=1}^{j=k}\frac{D}{\d t}(h_{ij}\ze_j)\Big|_{(b,t)}=0~~\forall\,i,t.$$
The second equation is equivalent to
\BE{ptrec_e3} 
\sum_{j=1}^{j=k}h_{ij}'(t)\ze_j(b,t)+
\sum_{j=1}^{j=k}\sum_{l=1}^{l=k}h_{ij}(t)A_{jl}(b,t)\ze_l(b,t)=0 
\qquad\Llra\qquad h'=-hA(b,\cdot).  \EE
Since (the real part of) the trace of $(A_{ij})$ is zero by \e_ref{ptrec_e1},
equation~\e_ref{ptrec_e3} has a unique solution in $\SO_k$ (or $\U_k$)  
such that $h(c)\!=\!\I$. 
Furthermore, by~\e_ref{ptrec_e2}
\BE{ptrec_e4}
\big|h(d)-\I\big|\le\int_c^d\!\!|h'(t)|\d t\le \int_c^d\!\!|h||A|\d t
\le k^2\!\int_c^d\!\!\!\int_a^b\!C_K|u_s||u_t|\d s\d t.\EE
Since $\Pi_{\partial\al}v_i\!=\!\sum_{j=1}^{j=k}h_{ij}(d)v_j$ by the above,
the lemma follows from equation~\e_ref{ptrec_e4}.

\begin{crl}\label{ptloop_crl}
If $(M,g)$ is a Riemannian manifold 
and $(E,\lr{,},\na)$ is a normed vector bundle with connection over~$M$,
for every compact subset $K\!\subset\!M$ there exists $C_K\!\in\!\R^+$ 
such that for every smooth closed curve $\al\!:[a,b]\!\lra\!M$ 
with $\Im\,\al\!\subset\!K$
$$\big|\Pi_{\al}-\I\big|\le C_K\min\big(\|\d\al\|_1,(b\!-\!a)\|\d\al\|_2^2\big).$$
\end{crl}

\noindent
{\it Proof:} Let $\exp\!:TM\!\lra\!M$ be an exponential-like map.
Since the group $\SO_k$ (or $\U_k$ if $E$ is complex) is compact and 
$$\|\d\al\|_1^2 \le (b\!-\!a)\|\d\al\|_2^2$$ 
by H\"older's inequality,
it is enough to assume that 
$$\|\d\al\|_1\le\min(r_{\exp}^g(K)/2,1).$$ 
Thus, there exists
$$\ti\al\in C^{\i}\big([a,b];T_{\al(a)}M\big) \qquad\hbox{s.t.}\qquad
\al(t)=\exp(\ti\al(t)), \quad |\ti\al(t)|_{\al(a)}<r_{\exp}(\al(a)).$$
Define 
$$u\!: [0,1]\!\times\![a,b]\!\lra\! K\subset M 
\qquad\hbox{by}\quad
u(s,t)=\exp\big(s\ti\al(t)\big).$$
Using
\begin{equation*}\begin{split}
|\ti\al(t)|&\le C_Kd_g\big(\al(a),\al(t)\big) \le C_K\|\d\al\|_1\,, \\
|\ti\al'(t)|&=\big|\{\d_{\ti\al(t)}\exp\}^{-1}(\al'(t))\big|
\le C_K|\d_t\al|\,,
\end{split}\end{equation*}
we find that 
\begin{alignat}{2}
\label{ptloop_e5a}
u_s(s,t)&=\big\{\d_{s\ti\al(t)}\exp\big\}\big(\ti\al(t)\big)
\qquad &\Lra& \qquad  |u_s|_{(s,t)}\le C_K'\|\d\al\|_1\,;\\ 
\label{ptloop_e5b}
u_t(s,t)&=s\big\{\d_{s\ti\al(t)}\exp\big\}\big(\ti\al'(t)\big)
\qquad &\Lra& \qquad |u_t|_{(s,t)}\le C_K'|\d_t\al|.
\end{alignat}
Thus, by Lemma~\ref{ptrec_lmm},
\begin{equation*}\begin{split}
\big|\Pi_{\al}-\I\big|
&=\big|\Pi_{\partial u}-\I\big|
\le  C_K\int_0^1\!\!\int_a^b\!\!|u_s||u_t|\d s\d t 
\le C_K'\|\d\al\|_1^2 \le C_K'(b\!-\!a)\|\d\al\|_2^2.
\end{split}\end{equation*}
Since $\|\d\al\|_1\!\le\!r_{\exp}^g(K)$, it follows that
$|\Pi_{\al}\!-\!\I|\!\le\!C_K\|\d\al\|_1$.

\begin{crl}\label{ptder_crl}
If $(M,g,\exp)$ is a Riemannian manifold with an exponential-like map
and $(E,\lr{,},\na)$ is a normed vector bundle with connection over~$M$,
for every compact subset $K\!\subset\!M$ there exists
 $C_K\!\in\!C^{\i}(\R;\R)$ such that for all $x\!\in\!K$ and smooth maps 
$\ti\al\!:(-\ep,\ep)\!\lra\!T_xM$  and $\xi\!:\!(-\ep,\ep)\!\lra\!E_x$
\BE{ptrec_e}\bigg|\frac{D}{\d t}\Big(\Pi_{\ti\al(t)}\xi(t)\Big)\Big|_{t=0}-
\Pi_{\ti\al(0)}\xi'(0)\bigg|\le 
C_K\big(|\ti\al(0)|\big)|\ti\al(0)||\ti\al'(0)||\xi(0)|.\EE
\end{crl}

\noindent
{\it Proof:} Define 
$$u\!:[0,1]\!\times\!\big[0,\ep/2\big] \lra K\subset M
\qquad\hbox{by}\quad  u(s,t)=\exp\big(s\ti\al(t)\big).$$ 
Let $\{v_i\}$ be an orthonormal basis for~$E_x$.
Extend each $v_i$ to 
$$\ze_i\in\Ga\big(u|_{[0,1]\times t};E\big)$$
by parallel-transporting along the curves $s\!\lra\!f(s,t)$. 
If 
$$\xi(t)=\sum_{i=1}^{i=k}f_i(t)v_i\,,$$ 
where $k$ is the rank of~$E$, then
\BE{ptder_e3}\begin{split}
\Pi_{\ti\al(t)}\xi(t)&=\sum_{i=1}^{i=k}f_i(t)\ze_i(1,t) \qquad\Lra\\
\frac{D}{\d t}\Big(\Pi_{\ti\al(t)}\xi(t)\Big)\Big|_{t=0}
&=\sum_{i=1}^{i=k}f_i'(0)\ze_i(1,0)+
\sum_{i=1}^{i=k}f_i(0)\frac{D}{\d t}\ze_i(1,t)\Big|_{t=0}\\
&=\Pi_{\ti\al(0)}\xi'(0)+\sum_{i=1}^{i=k}f_i(0)\frac{D}{\d t}\ze_i(1,t)\Big|_{t=0}.
\end{split}\EE
On the other hand, by \e_ref{ptrec_e0}, \e_ref{ptrec_e2},
and the first identities in~\e_ref{ptloop_e5a} and~\e_ref{ptloop_e5b},
\BE{ptder_e5}\begin{split}
\Big|\frac{D}{\d t}\ze_i(1,t)\Big|_{t=0}
&=\sum_{j=1}^{j=k}\big|A_{ij}(1,0)\big|
\le kC_K'\big(|\ti\al(0)|\big)\int_0^1\!\! |u_s|_{(s,0)}|u_t|_{(s,0)}\d s\\
&\le C_K\big(|\ti\al(0)|\big)|\ti\al(0)||\ti\al'(0)|.
\end{split}\EE
The claim follows from \e_ref{ptder_e3} and \e_ref{ptder_e5}.

\begin{rmk}\label{ptder_rmk}
Note that \e_ref{ptrec_e2} is applied with $K$ replaced by the compact set
$$\exp\big(\big\{v\!\in\!T_xM\!:\,x\!\in\!K,\,|v|\!\le\!|\ti\al(0)|\big\}\big);$$
thus, the constants $C_K'(|\ti\al(0)|)$ and $C_K(|\ti\al(0)|)$
 may depend on $|\ti\al(0)|$.
If $M$ is compact, then the first constant does not depend on $|\ti\al(0)|$,
since \e_ref{ptrec_e2} can then be applied with $K\!=\!M$.
The second constant is then also independent of $K$ and $|\ti\al(0)|$ 
if $\exp\!=\!\exp^{\na}$ for some connection~$\na$ in~$TM$.
So, in this case, the function $C_K$ in~\e_ref{ptrec_e} 
can be taken to be a constant independent of~$K$. 
\end{rmk}

\subsection{Poincare lemmas}
\label{poin_sect}

\begin{lmm}
\label{poin_lmm1}
If $\ze\!: S^1\!\lra\!\R^k$ is a smooth function such that
$\int_0^{2\pi}\!\ze(\th)\d\th\!=\!0$,
$$\int_0^{2\pi}\!\!|\ze(\th)|^2\d\th\le\int_0^{2\pi}\!\!|\ze'(\th)|^2\d\th.$$
\end{lmm}

\noindent
{\it Proof:} Write 
$$\ze(\th)=\sum_{n>-\i}^{n<\i}\!\ze_ne^{\fI n\th}\,;$$
see \cite[Section 6.16]{W}.
Since $\ze$ integrates to $0$, $\ze_0\!=\!0$. Thus,
$$\int_0^{2\pi}\!\!|\ze(\th)|^2\d\th = 2\pi\!\!\!\sum_{n>-\i}^{n<\i}\!|\ze_n|^2
\le 2\pi\!\!\!\sum_{n>-\i}^{n<\i}\!|n\ze_n|^2 =\int_0^{2\pi}\!\!|\ze'(\th)|^2\d\th.$$

\begin{prp}
\label{poincare_prp}
If $(M,g)$ is a Riemannian manifold and $(E,\lr{,},\na)$ is a normed vector bundle
with connection over~$M$,
for every compact subset $K\!\subset\!M$ there exists $C_K\!\in\!\R^+$ 
with the following property.
If $\al\!\in\!C^{\i}(S^1;M)$ is such that $\Im\,\al\!\subset\!K$
and $\xi,\ze\!\in\!\Ga(\al;E)$, then
$$\big|\llrr{\na_{\th}\xi,\ze}\big|\le
\|\na_{\th}\xi\|_2\|\na_{\th}\ze\|_2+
C_K\min\big(\|\d\al\|_1,\|\d\al\|_2^2\big)\|\xi\|_{2,1}\|\ze\|_2\,,$$
where $\na_{\th}\!\equiv\!\na^{\al}_{\partial_{\th}}$ is the covariant derivative
with respect to the oriented unit field on $S^1$ and
all the norms are computed with respect to the standard metric on~$S^1$.
\end{prp}

\noindent
{\it Proof:} Identify $E_{\al(0)}$ with $\R^k$ (or $\C^k$), 
preserving the metric.
Denote by $so(E_{\al(0)})\!\approx\!so_k$ (or $u(E_{\al(0)})\!\approx\!u_k$)
the Lie algebra of the Lie group $\SO(E_{\al(0)})\!\approx\!\SO_k$
(or of $\U(E_{\al(0)})\!\approx\!\U_k$).
For each $\chi\!\in\!so(E_{\al(0)})$ (or $\chi\!\in\!u(E_{\al(0)})$),
let $\E^{\chi}\!\in\!\SO(E_{\al(0)})$ (or $\E^{\chi}\!\in\!\U(E_{\al(0)})$)
be the exponential of~$\chi$.
Given $v\!\in\!E_{\al(0)}$, let $\ze_v(\th)\!\in\!E_{\al(\th)}$ denote 
the parallel transport of $v$ along the curve 
$t\!\lra\!\al(t)$ with $0\!\le\!t\!\le\!\th$.
By Corollary~\ref{ptloop_crl},
there exists $\chi\!\in\!so(E_{\al(0)})$ (or $\chi\!\in\!u(E_{\al(0)})$) such that 
\BE{poin_e0}|\chi|\le C_K\min\big(\|\d\al\|_1,\|\d\al\|_2^2\big)
\qquad\hbox{and}\qquad
\ze_v(2\pi)=\E^{\chi}\big(\ze_v(0)\big)=\E^{\chi}(v)
~~\forall\,v\!\in\!E_{\al(0)}\,.\EE
By the second statement in~\e_ref{poin_e0},
$$\Psi\!: S^1\!\times\!E_{\al(0)}\lra\al^*E\,, \qquad
(\th,v)\lra \ze_{\E^{-\th\chi/2\pi}(v)}(\th)\,,$$
is a smooth isometry.
Let $\Phi_2\!=\!\pi_2\!\circ\!\Psi^{-1}\!:\al^*E\lra E_{\al(0)}$ and
$$\bar\ze=\frac{1}{2\pi}\int_0^{2\pi}\{\Phi_2\ze\}(\th)\d\th\in E_{\al(0)}.$$
By H\"older's inequality and Lemma~\ref{poin_lmm1}, 
\BE{poin_e1}\begin{split}
\big|\llrr{\na_{\th}\xi,\ze\!-\!\Psi\bar\ze}\big|
&\le\|\na_{\th}\xi\|_2\|\ze\!-\!\Psi\bar\ze\|_2\\
&=\|\na_{\th}\xi\|_2\|\Phi_2\ze\!-\!\bar\ze\|_2
\le \|\na_{\th}\xi\|_2\|\d(\Phi_2\ze)\|_2.
\end{split}\EE
Note that
\BE{poin_e2}\begin{split}
\|\d(\Phi_2\ze)\|_2
&\le \|\na_{\th}\ze\|_2+|\chi/2\pi|\|\ze\|_2\\
&\le \|\na_{\th}\ze\|_2+C_K\min\big(\|\d\al\|_1,\|\d\al\|_2^2\big)\|\ze\|_2.
\end{split}\EE
On the other hand, by integration by parts, we obtain
\BE{poin_e3}
\llrr{\na_{\th}\xi,\ze\!-\!\Psi\bar\ze}
=\llrr{\na_{\th}\xi,\ze}+ \llrr{\xi,\na_{\th}(\Psi\bar\ze)}.  \EE
Since $\Psi\bar\ze$ is the parallel transport of 
$\E^{\th\chi/2\pi}\bar\ze$,
\BE{poin_e4}\begin{split}
\big|\llrr{\xi,\na_{\th}(\Psi\bar\ze)}\big|
&\le \|\xi\|_2\|\na_{\th}(\Psi\bar\ze)\|_2
=\|\xi\|_2|\chi/2\pi|\big\|\Psi\bar\ze\big\|_2\\
&\le C_K\min\big(\|\d\al\|_1,\|\d\al\|_2^2\big)\|\xi\|_2\|\ze\|_2.
\end{split}\EE
The proposition follows from equations \e_ref{poin_e1}-\e_ref{poin_e4}.\\

\noindent
Let $B_{R,r}\!\subset\!\R^2$ denote the open annulus with radii $r\!<\!R$
centered at the origin.

\begin{crl}[of Lemma~\ref{poin_lmm1}]\label{ellannbd_crl0}
There exists $C\!\in\!C^{\i}(\R;\R)$ such that for all $R\!\in\!\R^+$
$$r\!\in\!(0,R],~~~ \ze\!\in\!C^{\i}\big(B_{R,r};\R^k\big),~~~
\int_{B_{R,r}}\!\!\ze=0   \qquad\Lra\qquad
\|\ze\|_1\le C(R/r)R^2\|\d\ze\|_2.$$
\end{crl}

\noindent
{\it Proof:} It is sufficient to assume that $k\!=\!1$.
Define
$$\xi\!:S^1\lra\R \qquad\hbox{by}\qquad
\xi(\th)=\int_r^R\ze(\rho,\th)\rho\d\rho.$$
By H\"older's inequality and Lemma~\ref{poin_lmm1}, 
\BE{ellannbd0_e1}\begin{split}
\bigg(\int_0^{2\pi}\bigg|\int_r^R\ze(\rho,\th)\rho\d\rho\bigg|\d\th \bigg)^2
&\le  2\pi\int_0^{2\pi}\big|\xi(\th)\big|^2\d\th
\le 2\pi\int_0^{2\pi}\big|\xi'(\th)\big|^2\d\th\\
&\le 2\pi\int_0^{2\pi}\bigg(\int_r^R\big|\d_{(\rho,\th)}\ze\big|\rho^2\d\rho\bigg)^2\d\th\\
& \le \frac{\pi R^4}{2}
\int_0^{2\pi}\int_r^R\big|\d_{(\rho,\th)}\ze\big|^2\rho\d\rho\d\th
=\frac{\pi R^4}{2}\|\d\ze\|_2^2\,.
\end{split}\EE
If the function $\rho\!\lra\!\ze(\rho,\th)$ does not change sign on $(r,R)$, then
$$\int_r^R\big|\ze(\rho,\th)\big|\rho\d\rho
= \bigg|\int_r^R\ze(\rho,\th)\rho\d\rho\bigg|.$$
On the other hand, if this function vanishes somewhere on $(r,R)$, then
$$\big|\ze(\rho,\th)\big|\le \int_r^R\big|\d_{(t,\th)}\ze\big|\d t~~\forall\,\rho
\quad\Lra\quad
\int_r^R\big|\ze(\rho,\th)\big|\rho\d\rho
\le \frac{R^2}{2}\int_r^R\big|\d_{(t,\th)}\ze\big|\d t\,.$$
Combining these two cases and using~\e_ref{ellannbd0_e1}
and H\"older's inequality, we obtain
\BE{ellannbd0_e2}\begin{split}
\int_0^{2\pi}\!\!\!\int_r^R\big|\ze(\rho,\th)\big|\rho\d\rho\d\th
&\le \int_0^{2\pi}\bigg|\int_r^R\ze(\rho,\th)\rho\d\rho\bigg|\d\th
+\frac{R^2}{2}\int_0^{2\pi}\!\!\!\int_r^R\big|\d_{(\rho,\th)}\ze\big|\d\rho\d\th\\
& \le \frac{\sqrt{\pi} R^2}{\sqrt{2}}\|\d\ze\|_2+
\frac{R^2}{2}\|\d\ze\|_2
\bigg(\int_0^{2\pi}\!\!\!\int_r^R\rho^{-1}\d\rho\d\th\bigg)^{1/2}\\
&= \sqrt{\frac{\pi}{2}}\Big(1+\sqrt{\ln(R/r)}\Big)R^2\|\d\ze\|_2\,.
\end{split}\EE

\begin{rmk}\label{l1bnd_rmk}
By Corollary~\ref{ellannbd_crl} below, $C$ can in fact be chosen to
be a constant function.
Corollary~\ref{ellannbd_crl0} suffices for gluing $J$-holomorphic maps 
in symplectic topology, but Corollary~\ref{ellannbd_crl} leads to a sharper version 
of Proposition~\ref{elli_prp1}; see Remark~\ref{ellicrl_rmk}.
\end{rmk}

\subsection{Exponential-like maps and differentiation}
\label{expdiff_subs}

\noindent
Let $(M,g,\exp,\na)$ be a smooth Riemannian manifold with an exponential-like map~$\exp$
and connection~$\na$ in~$TM$, which is $g$-compatible, but not necessarily torsion-free.
Let
$$T_{\na}(\xi(x),\ze(x)\big)\equiv 
\big(\na_{\xi}\ze-\na_{\ze}\xi-[\xi,\ze]\big)\big|_x
\qquad\forall\,x\!\in\!M,\,\xi,\ze\!\in\!\Ga(M;TM),$$
be the torsion tensor of $\na$.
If $\al\!\!:\!(-\ep,\ep)\!\lra\!M$ is a smooth curve and $\xi\!\in\!\Ga(\al;TM)$, put
$$\Phi_{\al(0)}\Big(\al'(0);\xi(0),\frac{D}{\d s}\xi\Big|_{s=0}\Big)
=\Pi_{\xi(0)}^{-1}
\bigg(\frac{\d}{\d s}\exp\big(\xi(s)\big)\Big|_{s=0}\bigg)
=\Pi_{\xi(0)}^{-1}\big(\{\d_{\xi(0)}\exp\}(\xi'(0))\big),$$
where $\xi'(0)\!\in\!T_{\xi(0)}(TM)$ is the tangent vector to the curve
$\xi\!:(-\ep,\ep)\!\lra\!TM$ at $s\!=\!0$.

\begin{lmm}
\label{expdiff_lmm}
If $(M,g,\exp,\na)$ is a smooth Riemannian manifold with an exponential-like map
and a $g$-compatible connection, there exists $C\!\in\!C^{\i}(TM;\R)$ 
such that 
$$\Big|\Phi_x(v;w_0,w_1)-\big(v\!+\!w_1\!-\!T_{\na}(v,w_0)\big)\Big|
\le C(w_0)\big(|v||w_0|^2\!+\!|w_0||w_1|\big)$$
for all $x\!\in\!M$ and $v,w_0,w_1\!\in\!T_xM$.
\end{lmm}

\noindent
{\it Proof:} Let $\al\!:(-\ep,\ep)\!\lra\!M$ be a smooth curve and 
$\xi\!\in\!\Ga(\al;TM)$ such~that
$$\al(0)=x,\quad\al'(0)=v,\quad
\xi(0)=w_0,\quad\frac{D}{\d s}\xi(s)\Big|_{s=0}=w_1.$$
Put
\begin{equation*}\begin{split}
F_{v,w_0,w_1}(t)&=\frac{\d}{\d s}\exp\big(t\xi(s)\big)\Big|_{s=0}
=\{\d_{tw_0}\exp\}\big(\d_{w_0}m_t(\xi'(0))\big), \\
H_{v,w_0,w_1}(t)&=\Pi_{tw_0}\big(v\!+\!tw_1\!-\!tT_{\na}\big(v,w_0)\big),
\end{split}\end{equation*}
where $m_t\!:TM\!\lra\!TM$ is the scalar multiplication by $t$.
Then,
\begin{gather*}
F_{v,w_0,w_1}(0)=\frac{\d}{\d s}\al(s)\Big|_{s=0}=v=H_{v,w_0,w_1}(0),\\
\frac{D}{\d t}F_{v,w_0,w_1}(t)\Big|_{t=0}=
\frac{D}{\d s}\frac{\d}{\d t}\exp\big(t\xi(s)\big)\Big|_{t=0}\Big|_{s=0}
-T_{\na}\big(v,w_0\big)
=w_1-T_{\na}(v,w_0)=\frac{D}{\d t}H_{v,w_0,w_1}(t)\Big|_{t=0};
\end{gather*}
see Corollary~\ref{ptder_crl}.  
Since 
$$F_{\cdot,w_0,\cdot}(t)-H_{\cdot,w_0,\cdot}(t)
\in\Hom(T_xM\!\oplus\!T_xM,T_{\exp(tw_0)}M),$$
combining the last two equations, we obtain
$$\big|F_{v,w_0,w_1}(t)-H_{v,w_0,w_1}(t)\big|\le C(w_0,t)t^2\big(|v|\!+\!|w_1|\big)
\quad\forall~v,w_0,w_1\!\in\!T_xM,~x\!\in\!M,~t\!\in\!\R,$$
where $C$ is a smooth function on $TM\!\times\!\R$.
Since 
$$F_{v,w_0,w_1}(t)-H_{v,w_0,w_1}(t)=F_{v,tw_0,tw_1}(1)-H_{v,tw_0,tw_1}(1),$$
we conclude that there exists $C\!\in\!C^{\i}(TM)$ such~that
\begin{equation}\label{diff1_e2}
\big|F_{v,w_0,w_1}(1)-H_{v,w_0,w_1}(1)\big|\le C(w_0)\big(|w_0|^2|v|\!+\!|w_0||w_1|\big)
\quad\forall~v,w_0,w_1\!\in\!T_xM,\,x\!\in\!M,
\end{equation}
as claimed.\\

\noindent 
For any $v,w_0,w_1\!\in\!T_xM$, let 
$\ti\Phi_x(v;w_0,w_1)=\Phi_x(v;w_0,w_1)-\big(v\!+\!w_1\!-\!T_{\na}(v,w_0)\big)$.

\begin{crl}\label{expdiff_crl}
If $(M,g,\exp,\na)$ is a smooth Riemannian manifold with an exponential-like map
and a $g$-compatible connection, there exists $C\!\in\!C^{\i}(TM\!\times_M\!TM;\R)$ 
such~that 
\begin{equation*}\begin{split}
&\Big|\ti\Phi_x(v;w_0,w_1)\!-\!\ti\Phi_x(v;w_0',w_1')\Big|\\
&\qquad\qquad\qquad\le
C(w_0,w_0')\Big(\! \big( (|w_0|+|w_0'|)|v|\!+\!|w_1|\!+\!|w_1'|\big)|w_0\!-\!w_0'|
+\big(|w_0|\!+\!|w_0'|\big)|w_1\!-\!w_1'|\Big)
\end{split}\end{equation*}
for all $x\!\in\!M$ and $v,w_0,w_1,w_0',w_1'\!\in\!T_xM$.
\end{crl}

\noindent
{\it Proof:} By the proof of Lemma~\ref{expdiff_lmm}, 
$$\ti\Phi(v;w_0,w_1)=\ti\Phi_1(w_0;v)+\ti\Phi_2(w_0;w_1)$$
for some smooth bundle sections
$\ti\Phi_1,\ti\Phi_2\!: TM\lra \pi_{TM}^*\Hom(TM,TM)$ such that
$$\big|\ti\Phi_1(w_0;\cdot)\big|\le C_1(w_0)|w_0|^2\,,\qquad
\big|\ti\Phi_2(w_0;\cdot)\big|\le C_2(w_0)|w_0| \qquad\forall~w_0\!\in\!TM.$$
Thus,
\begin{equation*}\begin{split}
\big|\ti\Phi_1(w_0;\cdot)-\ti\Phi_1(w_0';\cdot)\big|
&\le C_1'(w_0,w_0')\big(|w_0|\!+\!|w_0'|\big)|w_0\!-\!w_0'|\\
\big|\ti\Phi_2(w_0;\cdot)-\ti\Phi_2(w_0';\cdot)\big|
&\le C_2'(w_0,w_0')|w_0\!-\!w_0'|
\end{split}
\qquad\forall~w_0,w_0'\!\in\!T_xM.
\end{equation*}
From the linearity of $\ti\Phi_1(w_0;\cdot)$ and $\ti\Phi_2(w_0;\cdot)$
in the second input, we conclude that
\begin{equation*}\begin{split}
\big|\ti\Phi_1(w_0;v)-\ti\Phi_1(w_0';v)\big|
&\le C_1'(w_0,w_0')\big(|w_0|\!+\!|w_0'|\big)|w_0\!-\!w_0'||v|,\\
\Big|\ti\Phi_2(w_0;w_1)-\ti\Phi_2(w_0;w_1')\Big|
&\le C_2'(w_0,w_0')|w_0\!-\!w_0'||w_1|+C_2(w_0')|w_0'||w_1\!-\!w_1'|.
\end{split}\end{equation*}

\subsection{Expansion of the $\dbar$-operator}

\noindent
Let $(M,J)$ and  $(\Si,\fj)$  be almost-complex manifolds.
If $u\!:\Si\!\lra\!M$ is a smooth map, let 
\begin{gather*}
\Ga(u)=\Ga(\Si;u^*TM),\qquad
\Ga_{J,\fj}^{0,1}(u)=\Ga\big(\Si;T^*\Si^{0,1}\!\otimes\!_{\C}u^*TM\big),\\
\dbar_{J,\fj}u=\frac{1}{2}\big(\d u+J\circ\d u\circ\fj\big)\in \Ga_{J,\fj}^{0,1}(u),
\end{gather*}
as in~\e_ref{dbarJjdfn_e}.
If $\na$ is a connection in $TM$, define
$$D_{J,\fj;u}^{\na}\!:\Ga(u)\lra \Ga_{J,\fj}^{0,1}(u)
\qquad\hbox{by}\qquad
D_{J,\fj;u}^{\na}\xi=\frac{1}{2}\big(\na^u\xi+J\na_{\fj}^u\xi\big)-
\frac{1}{2}\big(T_{\na}(\d u,\xi)+JT_{\na}(\d u\!\circ\!\fj,\xi)\big).$$
If in addition $\exp\!:TM\!\lra\!M$ is an exponential-like map
and $\na J\!=\!0$, define
\begin{gather*}
\exp_u\!:\Ga(u)\lra C^{\i}(\Si;M),\quad
\dbar_u,N_{\exp}^{\na}\!:\Ga(u)\lra \Ga_{J,\fj}^{0,1}(u) \qquad\hbox{by}\\
\big\{\!\exp_u(\xi)\big\}(z)=\exp\!\big(\xi(z)\big)~~\forall\,z\!\in\!\Si,
\quad
\big\{\dbar_u\xi\big\}_z(v)=
\Pi_{\xi(z)}^{-1}\big(\big\{\dbar_{J,\fj}(\exp_u(\xi))\big\}_z(v)\big)
~~\forall\,z\!\in\!\Si,\,v\!\in\!T_z\Si,\\\
\dbar_u\xi=\dbar_{J,\fj}u+D_{J,\fj;u}^{\na}\xi+N_{\exp}^{\na}(\xi).
\end{gather*}

\begin{lmm}\label{dbar_lmm}
If $(M,J,g,\exp,\na)$ is an almost-complex Riemannian manifold with 
an exponential-like map and a $g$-compatible connection in~$(TM,J)$, 
there exists $C\!\in\!C^{\i}(TM\!\times_{M}\!TM;\R)$ with the following property.
If $(\Si,\fj)$ is an almost complex manifold, $u\!:\Si\!\lra\!M$ is a smooth map,
and $\xi,\xi'\!\in\!\Ga(u)$, then 
\begin{equation*}\begin{split}
\Big|\big\{N_{\exp}^{\na}(\xi)\big\}_z(v)-\big\{N_{\exp}^{\na}(\xi')\big\}_z(v)\Big|
\le C\big(\xi(z),\xi'(z)\big) \Big(
\big(|\xi(z)|\!+\!|\xi'(z)|\big)
\big(|\na_v(\xi\!-\xi')|+|\na_{\fj v}(\xi\!-\xi')|\big)&\\
+\big((|\d_zu(v)|\!+\!|\d_zu(\fj v)|)(|\xi(z)|\!+\!|\xi'(z)|)
+(|\na_v\xi|\!+\!|\na_{\fj v}\xi|+|\na_v\xi'|\!+\!|\na_{\fj v}\xi|)\big)
\big|\xi(z)\!-\!\xi'(z)\big|\Big)&
\end{split}\end{equation*}
for all  $z\!\in\!\Si$, $v\!\in\!T_z\Si$.
Furthermore, $N_{\exp}^{\na}(0)\!=\!0$.
\end{lmm}

\noindent
{\it Proof:}
Since the connection $\na$ commutes with $J$, 
so does the parallel transport~$\Pi$.
Thus, with notation as in Section~\ref{expdiff_subs},
$$\big\{N_{\exp}^{\na}(\xi)\big\}_z(v)=
\frac{1}{2}\Big(
\ti\Phi\big(\d_zu(v);\xi(z),\na_v\xi\big)
+J\big(u(z)\big)\ti\Phi\big(\d_zu(\fj v);\xi(z),\na_{\fj v}\xi\big)\Big).$$
The claim now follows from Corollary~\ref{expdiff_crl}.

\begin{dfn}\label{norms_dfn}
Let $M$ be a smooth manifold and 
$(E,\lr{,},\na)$ a normed vector bundle with connection over~$M$.
If $C_0\!\in\!\R^+$, $(\Si,\fj)$ is an almost complex manifold,
and $u\!: \Si\!\lra\!M$ is a smooth map, 
norms $\|\cdot\|_{p,1}$ and $\|\cdot\|_p$  on $\Ga(u;E)$ and $\Ga^1(u;E)$, 
respectively, are \sf{$C_0$-admissible} if for all 
$\xi\!\in\!\Ga(u;E)$, $\eta\!\in\!\Ga^1(u;E)$,
and every continuous function $f\!:\Si\!\lra\!\R$,
$$\|f\eta\|_p\le\|f\|_{C^0}\|\eta\|_p,\quad
\|\eta\circ\fj\|_p=\|\eta\|_p,\quad
\|\na^u\xi\|_p\le\|\xi\|_{p,1},\quad
\|\xi\|_{C^0}\le C_0\|\xi\|_{p,1}.$$
\end{dfn}

\begin{prp}\label{dbar_prp}
If $(M,J,g,\exp,\na)$ is an almost-complex Riemannian manifold with 
an\linebreak 
exponential-like map and a $g$-compatible connection in~$(TM,J)$, 
for every compact subset $K\!\subset\!M$ 
there exists $C_K\!\in\!C^{\i}(\R;\R)$ with the following property.
If $(\Si,\fj)$ is an almost complex manifold, $u\!:\Si\!\lra\!M$ is a smooth map,
and $\|\cdot\|_{p,1}$ and $\|\cdot\|_p$ are 
$C_0$-admissible norms on $\Ga(u;TM)$ and $\Ga^1(u;TM)$, respectively,
then 
$$\big\|N_{\exp}^{\na}(\xi)-N_{\exp}^{\na}(\xi')\big\|_p
\le C_K\big(C_0\!+\!\|\d u\|_p\!+\!\|\xi\|_{p,1}\!+\!\|\xi'\|_{p,1}\big)
\big(\|\xi\|_{p,1}\!+\!\|\xi'\|_{p,1}\big)\|\xi\!-\!\xi'\|_{p,1}$$
for all $\xi,\xi'\!\in\!\Ga(u)$.
Furthermore, $N_{\exp}^{\na}(0)\!=\!0$.
If the $g$-ball $B_{g;\de}(u(z))$ of radius~$\de$ around $f(z)$ for some $z\!\in\!\Si$
is isomorphic to an open subset of $\C^n$ and $|\xi(z)|\!<\!\de$, 
then $\{N_{\exp}^{\na}\xi\}_z\!=\!0$.
\end{prp}

\noindent
{\it Proof:} The first two statements follow
from Lemma~\ref{dbar_lmm} and Definition~\ref{norms_dfn}.
The last claim is clear from the definition of~$N_{\exp}^{\na}$.

\begin{rmk}\label{admissnorms_rmk}
As the notation suggests, 
one possibility for the norms $\|\cdot\|_{p,1}$ and $\|\cdot\|_p$ is
the usual Sobolev $L^p_1$ and $L^p$-norms with respect
to some Riemannian metric on~$\Si$, where $p\!>\!\dim_{\R}\Si$.
Another natural possibility in the $\dim_{\R}\Si\!=\!2$ case is 
the modified Sobolev norms introduced in \cite[Section~3]{LT};
these are particularly suited for gluing pseudo-holomorphic curves.
By Proposition~\ref{c0bound_prp} below,  in the $\dim_{\R}\Si\!=\!2$ case 
the constant $C_0$ itself is a function of $\|\d u\|_p$~only
for either of these two choices of norms.
\end{rmk}

\begin{rmk}\label{Nijen_rmk}
By Proposition~\ref{dbar_prp}, the operator $D_{J,\fj;u}^{\na}$ defined above
is a linearization of the $\dbar$-operator on the space of smooth maps to~$M$ at~$u$.
If $\na'$ is any connection in $TM$, the connection
$$\na\!:\Ga(M;TM)\lra\Ga(M;T^*M\!\otimes_{\R}\!TM),\quad
\na_v\xi=\frac{1}{2}\Big(\na_v'\xi-J\na_v'(J\xi)\Big)
\quad\forall\,v\!\in\!TM,\,\xi\!\in\!\Ga(M;TM),$$
is $J$-compatible.
If in addition $\na'$ and $J$ are compatible with a Riemannian metric $g$ on $M$,
then so is~$\na$.
If $\na'$ is also the Levi-Civita connection of the metric~$g$ (i.e.~$T_{\na'}\!=\!0$),
$$T_{\na}(v,w)=\frac{1}{2}\big(J(\na_w'J)v-J(\na_v'J)w\big)
\qquad \forall\,v,w\!\in\!T_xM,\,x\!\in\!M.$$
If the two-form $\om(\cdot,\cdot)\!\equiv\!g(J\cdot,\cdot)$ is  closed as well,
then 
$$\na_{Jv}'J=-J\na_v'J \qquad\forall\,v\!\in\!TM$$
by \cite[(C.7.5)]{MS} and thus 
$$T_{\na}(v,w)=
-\frac{1}{4}\big(J(\na_v'J)w-J(\na_w'J)v-(\na_{Jv}'J)w+(\na_{Jw}'J)v\big)
= -A_J(v,w) \quad \forall\,v,w\!\in\!T_xM,\,x\!\in\!M,$$
where $A_J$ is the Nijenhuis tensor of $J$ as in~\e_ref{Nijendfn_e}.
The operator $D_{J,\fj;u}^{\na}$ then becomes
\BE{Dcomm_e}
D_{J,\fj;u}^{\na}\!:\Ga(u)\lra \Ga_{J,\fj}^{0,1}(u),  \qquad
D_{J,\fj;u}^{\na}\xi=\dbar_{\na^u}\xi+A_J(\partial_{J,\fj}u,\xi),\EE
where
\begin{equation*}\begin{split}
\dbar_{\na^u}\xi&=\frac{1}{2}\big(\na^u\xi+J\na_{\fj}^u\xi\big)
\in\Ga_{J,\fj}^{0,1}(u),\\
\partial_{J,\fj}u&=\frac{1}{2}\big(\d u-J\circ \d u\circ\fj\big)
\in\Ga\big(\Si;T^*\Si^{1,0}\!\otimes_{\C}\!u^*TM\big).
\end{split}\end{equation*}
This agrees with \cite[(3.1.5)]{MS}, since the Nijenhuis tensor of $J$
is defined to be $-4A_J$ in \cite[p18]{MS}.
\end{rmk}

\section{Sobolev and elliptic inequalities}
\label{sob_sec}

\noindent
This appendix refines, in the $n\!=\!2$ case, 
the proofs of Sobolev Embedding Theorems given in~\cite{M}
to obtain a $C^0$-estimate in Proposition~\ref{c0bound_prp}
and elliptic estimates for the $\dbar$-operator
in Propositions~\ref{elli_prp1} and~\ref{elli_prp2}.
If $R,r\!\in\!\R$, let
$$B_R=\big\{x\!\in\!\R^2\!:\,|x|\!<\!R\big\},\qquad
B_{R,r}=B_R-\bar{B}_r\,,\qquad \ti{B}_{R,r}=B_R-B_r\,.$$

\subsection{Eucledian case}

\noindent
If $\xi$ is an $\R^k$-valued function defined on a subset $B$ of $\R^2$,
let $\supp_{\R^2}(\xi)$ be the closure of $\supp(\xi)\!\subset\!B$ in~$\R^2$.
If $U$ is an open subset of $\R^2$, $\xi\!\in\!C^{\i}(U;\R^k)$, 
and $p\!\ge\!1$, let
$$\|\xi\|_p\equiv\bigg(\int_U|\xi|^p\bigg)^{1/p}\,,\qquad
\|\xi\|_{p,1}\equiv\|\xi\|_p+\|\d\xi\|_p\,,$$
be the usual Sobolev norms of $\xi$.

\begin{lmm}\label{c0plane_lmm}
For every bounded convex domain $\cD\!\subset\!\R^2$,
$\xi\!\in\!C^{\i}(\cD;\R^k)$, and $x\!\in\!\cD$,
$$\big|\xi_{\cD}-\xi(x)\big|\le 
\frac{2r_0^2}{|\cD|} \int_{\cD}\!|\d_y\xi||y\!-\!x|^{-1}\d y,$$
where $2r_0$ is the diameter of $\cD$, 
$|\cD|$ is the area of $\cD$, and 
$$\xi_{\cD}=\frac{1}{|\cD|}\Big(\int_{\cD}\xi(y)\d y\Big)$$
is the average value of $\xi$ on $\cD$.
\end{lmm}

\noindent
{\it Proof:} For any $y\!\in\!\cD$,
$$\xi(y)-\xi(x)=\int_0^1\!\frac{\d}{\d t}\xi\big(x\!+\!t(y\!-\!x)\big)\d t
=\int_0^1\! \d_{x+t(y-x)}\xi(y\!-\!x) \d t.$$
Putting $g(z)\!=\!|\d_z\xi|$ if $z\!\in\!\cD$ and $g(z)\!=\!0$ otherwise, we obtain 
$$\big|\xi_{\cD}-\xi(x)\big|\le \frac{1}{|\cD|}
\int_{y\in\cD}|\xi(y)\!-\!\xi(x)|\d y
\le\frac{1}{|\cD|}\int_{y\in\cD}\int_0^{\i}
g\big(x\!+\!t(y\!-\!x)\big)|y\!-\!x|\d t\d y.$$
Rewriting the last integral in polar coordinates $(r,\th)$ centered at $x$, 
we obtain
\begin{equation*}\begin{split}
\big|\xi_{\cD}-\xi(x)\big| &\le\frac{1}{|\cD|}
\int_0^{2\pi}\!\!\!\int_0^{2r_0}\!\!\!\int_0^{\i}\!\!g(tr,\th)r^2\d t\d r\d\th\\
&=\frac{1}{|\cD|}\int_0^{2\pi}\!\!\!\int_0^{2r_0}\!\!\!\int_0^{\i}\!\!
g(t,\th)r\d t\d r\d\th=
\frac{2r_0^2}{|\cD|}\int_0^{2\pi}\!\!\!\int_0^{\i}\!\!g(t,\th)\d t\d\th\\
&=\frac{2r_0^2}{|\cD|}\int_{\cD}|\d_y\xi||y\!-\!x|^{-1}\d y.
\end{split}\end{equation*}

\begin{crl}\label{c0plane_crl0}
For every $p\!>\!2$, there exists $C_p\!>\!0$ such~that
$$r\!\in\!\big[0,R/2\big],~~\xi\!\in\!C^{\i}(B_{R,r};\R^k) 
\qquad\Lra\qquad 
\big|\xi(x)-\xi(y)\big|\le C_pR^{\frac{p-2}{p}}\|\d\xi\|_p
~~~\forall\,x,y\!\in\!B_{R,r}\,.$$
\end{crl}

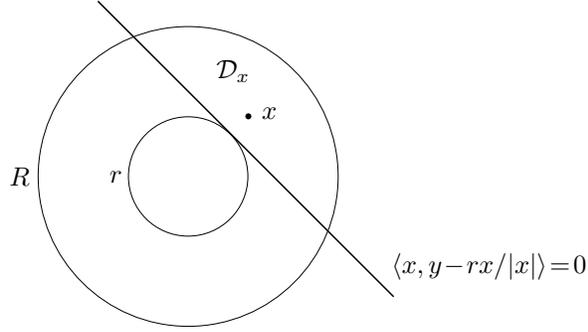
\begin{figure}
\begin{pspicture}(0,-2)(10,2.5)
\psset{unit=.4cm}
\pscircle[linewidth=.03](20,0){5}\rput(14.4,0){\sm{$R$}}
\pscircle[linewidth=.03](20,0){2}\rput(17.6,0){\sm{$r$}}
\pscircle*(22,2){.1}\rput(22.7,2.1){\sm{$x$}}
\psline[linewidth=.05](26.83,-4)(17,5.83)
\rput(30,-3){\sm{$\lr{x,y\!-\!rx/|x|}\!=\!0$}}
\rput(21.5,3.5){\sm{$\cD_x$}}
\end{pspicture}
\caption{A convex region $\cD_x$ of the annulus $\cD_{R,r}$ containing $x$} 
\label{ann_fig}
\end{figure}

\noindent
{\it Proof:} For any $x\!\in\!B_{R,r}$, put
$$\cD_x=\big\{y\!\in\!B_{R,r}\!: 
\lr{x,|x|y\!-\!rx}\!>\!0\big\}.$$
If $x\!\neq\!0$, $\cD_x$ is the part of the annulus on the same side 
of the line $\lr{x,y\!-\!rx/|x|}\!=\!0$ as~$x$;
see Figure~\ref{ann_fig}.
In particular,
$$\diam(\cD_x)\le 2R\,, \qquad
|\cD_x|\ge\Big(\frac{\pi}{3}\!-\!\frac{\sqrt{3}}{4}\Big)R^2.$$
Thus, by Lemma~\ref{c0plane_lmm} and H\"older's inequality,
\BE{C0diff_e1}\begin{split}
\big|\xi(x)-\xi_{\cD_x}|&\le 12\int_{y\in\cD_x}\!\!|\d_y\xi||y\!-\!x|^{-1}\d y\\
&\le 12\bigg(\int_{y\in B_{2R}(x)}\!|y\!-\!x|^{-\frac{p}{p-1}}\bigg)^{\frac{p-1}{p}}
\|\d\xi\|_p \le C_pR^{\frac{p-2}{p}}\|\d\xi\|_p,
\end{split}\EE
since $\frac{p}{p-1}\!<\!2$.
Let 
$$x_{\pm}=\big(\pm(R\!-\!r)/2,0\big), \qquad 
y_{\pm}\!=\!\big(0,\pm(R\!-\!r)/2\big).$$
Since each of the convex regions $\cD_{x_{\pm}}$ intersects 
$\cD_{y_+}$ and $\cD_{y_-}$ and 
$\cD_x$ intersects at least one (in fact precisely two if $r\!\neq\!0$) 
of these four convex regions for every $x\!\in\!B_{R,r}$,
$$\big|\xi(x)-\xi(y)\big|\le 8C_pR^{\frac{p-2}{p}}\|\d\xi\|_p
~~~\forall\,x,y\!\in\!B_{R,r}$$
by~\e_ref{C0diff_e1} and triangle inequality.

\begin{crl}\label{c0plane_crl}
For every $p\!>\!2$, there exists $C_p\!\in\!C^{\i}(\R^+;\R)$ 
such~that
$$r\!\in\!\big[0,R/2\big],~~\xi\!\in\!C^{\i}(B_{R,r};\R^k) 
\quad\Lra\quad \|\xi\|_{C^0}\le C_p(R)\|\xi\|_{p,1}.$$
\end{crl}

\noindent
{\it Proof:} 
By Corollary~\ref{c0plane_crl0} and H\"older's inequality, for every $x\!\in\!B_{R,r}$
\BE{c0plane_e}\begin{split}
|\xi(x)|&\le \big|\xi_{B_{R,r}}\big|+C_pR^{\frac{p-2}{p}}\|\d\xi\|_p
\le \frac{1}{|B_{R,r}|}\|\xi\|_1+C_pR^{\frac{p-2}{p}}\|\d\xi\|_p\\
&\le |B_{R,r}|^{-\frac{1}{p}}\|\xi\|_p+C_pR^{\frac{p-2}{p}}\|\d\xi\|_p
\le(1\!+\!C_p)R^{-\frac{2}{p}}\big(\|\xi\|_p+R\|\d\xi\|_p\big).
\end{split}\EE

\begin{lmm}
\label{loj_lmm}
For all $R\!>\!0$ and $r\!\in\![0,R)$,
$$\ze\!\in\! C^{\i}(B_{R,r};\R^k),\quad
\supp_{\R^2}(\ze)\!\subset\!\ti{B}_{R,r} 
\quad\Lra\quad \|\ze\|_2\le\|\d\ze\|_1.$$
\end{lmm}

\noindent
{\it Proof:} Such a function $\ze$ can be viewed as a function on 
the complement of the ball $B_r$ in~$\R^2$.
Since $\ze$ vanishes at infinity, for any $(x,y)\!\in\!B_{R,r}$
$$\ze(x,y)=\begin{cases}
\int_{-\i}^x\ze_s(s,y)\d s,&\hbox{if~}x\!\le\! 0;\\
-\int_x^{\i}\ze_s(s,y)\d s,&\hbox{if~}x\!\ge\! 0;
\end{cases}
\qquad
\ze(x,y)=\begin{cases}
\int_{-\i}^y\ze_t(x,t)\d t,&\hbox{if~}y\!\le\! 0;\\
-\int_y^{\i}\ze_t(x,t)\d t,&\hbox{if~}y\!\ge\! 0.
\end{cases}$$
Taking the absolute value in these equations, we obtain
\begin{equation}
\label{loj_e1}
\big|\ze(x,y)\big|\le \int_{-\i}^{\i}\!\big|\d_{(s,y)}\ze\big|\d s
\quad\hbox{and}\quad
\big|\ze(x,y)\big|\le \int_{-\i}^{\i}\!\big|\d_{(x,t)}\ze\big|\d t,
\end{equation}
where we formally set $\ze$ and $\d\ze$ to be zero on the smaller disk.
Multiplying the two inequalities in \e_ref{loj_e1} and integrating with 
respect to $x$ and $y$, we conclude
$$\int_{-\i}^{\i}\!\int_{-\i}^{\i}\!\big|\ze(x,y)\big|^2\d x\d y\le
\Big(\int_{-\i}^{\i}\!\int_{-\i}^{\i}\!\big|\d_{(x,y)}\ze\big|\d x\d y\Big)^2,$$
as claimed.

\begin{crl}\label{pqplane_crl}
For all $p,q\!\ge\! 1$ with $1\!-\!2/p\ge -2/q$,
there exists $C_{p,q}\!\in\!\R^+$ such that
$$r\!\in\![0,R),\quad
\xi\!\in\! C^{\i}(B_{R,r};\R^k),\quad
\supp_{\R^2}(\xi)\!\subset\!\ti{B}_{R,r}
\quad\Lra\quad \|\xi\|_q\le C_{p,q}R^{1-\frac{2}{p}+\frac{2}{q}}\|\d\xi\|_p.$$
\end{crl}

\noindent
{\it Proof:} We can assume that $k\!=\!1$. For $\ep\!>\!0$, let 
$\ze_{\ep}=(\xi^2\!+\!\ep)^{\frac{q}{4}}-\ep^{\frac{q}{4}}$. 
By Lemma~\ref{loj_lmm} and H\"older's inequality,
\BE{pqplanebound_e1}\begin{split}
\|\xi\|_q^q &\le \big\|\ze_{\ep}\!+\!\ep^{\frac{q}{4}}\big\|_2^2\le 
2\|\d\ze_{\ep}\|_1^2+2\ep^{\frac{q}{2}}\pi R^2
=2\big\|\frac{q}{2}(\xi^2\!+\!\ep)^{\frac{q}{4}-1}\xi\d\xi\big\|_1^2+
  2\ep^{\frac{q}{2}}\pi R^2\\
&\le q^2\big\|(\xi^2\!+\!\ep)^{\frac{q}{4}-\frac{1}{2}}\d\xi\big\|_1^2
      +2\ep^{\frac{q}{2}}\pi R^2
\le q^2\|\d\xi\|_p^2\big\|(\xi^2\!+\!\ep)^{\frac{q-2}{4}}\big\|_{\frac{p}{p-1}}^2+
2\ep^{\frac{q}{2}}\pi R^2.
\end{split}\EE
Note that
$$1-\frac{2}{p}=-\frac{2}{q} \quad\Lra\quad 
\frac{q-2}{4}\frac{p}{p-1}=\frac{q-2}{4}\frac{2q}{q-2}=\frac{q}{2}.$$
Thus, letting $\ep$ go to zero in~\e_ref{pqplanebound_e1}, we obtain
$$\|\xi\|_q^q\le q^2 \|\d\xi\|_p^2\|\xi\|_q^{q-2}
\quad\Lra\quad \|\xi\|_q\le q\|\d\xi\|_p.$$
The case $1\!-\!\frac{2}{p}>-\frac{2}{q}$ follows by H\"older's inequality.

\begin{rmk}\label{pqplanecrl_rmk}
By H\"older's inequality, the constant $C_{p,q}$ can be taken to be
$$C_{p,q}=\max(2,q)\pi^{\frac{1}{2}\left(1-\frac{2}{p}+\frac{2}{q}\right)}\,.$$
\end{rmk}

\begin{crl}[of Lemmas~\ref{c0plane_lmm}, \ref{loj_lmm}]\label{ellannbd_crl}
There exists $C\!>\!0$ such that for all $R\!\in\!\R^+$
$$r\!\in\![0,R],~~~ \ze\!\in\!C^{\i}\big(B_{R,r};\R^k\big),~~~
\int_{B_{R,r}}\!\!\ze=0   \qquad\Lra\qquad
\|\ze\|_1\le CR^2\|\d\ze\|_2.$$
\end{crl}

\begin{figure}
\begin{pspicture}(0,-2)(10,2.5)
\psset{unit=.4cm}
\pscircle[linewidth=.03](20,0){6}\rput(12.7,0){\sm{$R\!=\!1$}}
\pscircle[linewidth=.03](20,0){1}\rput(18.6,0){\sm{$r$}}
\psline[linewidth=.05](19.13,.5)(14.8,3)
\psline[linewidth=.05](20.87,.5)(25.2,3)
\psline[linewidth=.05](20,-1)(20,-6)
\psline[linewidth=.05](18.27,1)(21.73,1)
\rput(20.2,3.5){\sm{$\cD_+$}}
\end{pspicture}
\caption{A large convex region $\cD_+$ of an annulus $\cD$} 
\label{ann_fig2}
\end{figure}
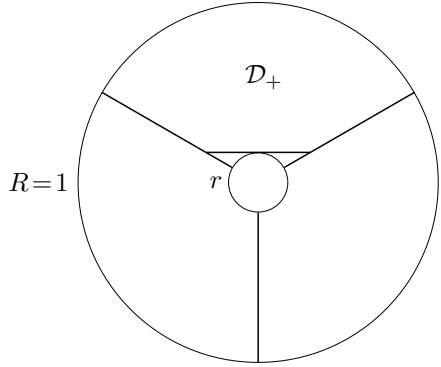

\noindent
{\it Proof:} (1) If $\ze\!\in\!C^{\i}(B_{R,r};\R^k)$ integrates to 0
over its domain, then so does the function
$$\ti\ze\!\in\!C^{\i}\big(B_{1,r/R};\R^k\big), \qquad
\ti\ze(z)=\ze(Rz).$$
Furthermore, $\|\ti\ze\|_1\!=\!\|\ze\|_1/R^2$ and 
$\|\d\ti\ze\|_2\!=\!\|\d\ze\|_2$.
Thus, it is sufficient to prove the claim for $R\!=\!1$.\\
(2) If $r\!=\!0$, for some open half-disk $\cD\!\subset\!B_{1,0}$
\BE{ellannbd_e1}\int_{\cD}\!\!\ze=0, \qquad 
\big\|\ze|_{\cD}\big\|_1\ge \frac{1}{2}\|\ze\|_1\,.\EE
By the first condition, Lemma~\ref{c0plane_lmm}, and H\"older's inequality
\begin{equation*}\begin{split}
\big\|\ze|_{\cD}\big\|_1\le
\frac{4}{\pi}\int_{\cD}\!\int_{\cD}|\d_y\ze||y\!-\!x|^{-1}\d y\d x
\le 16\int_{\cD}|\d_y\ze|\d y
\le 8\sqrt{2\pi}\|\d\ze\|_2\,.
\end{split}\end{equation*}
Along with the second assumption in~\e_ref{ellannbd_e1}, this  
implies the claim for $r\!=\!0$ with $C\!=\!16\sqrt{2\pi}$.\\
(3) Let $\be\!:\R\!\lra\![0,1]$ be a smooth function such that 
$$\be(t)=\begin{cases}1,&\hbox{if}~t\le1/2; \\
0,&\hbox{if}~t\ge1.\end{cases}$$
It remains to prove the claim for all $r\!>\!0$ and $R\!=\!1$.
By~\e_ref{ellannbd0_e2}, we can assume that 
\BE{ellannbd_e3} r\le\frac{1}{48\sqrt{3\pi}\|\be'\|_{C^0}}
<\frac{1}{96{\sqrt{3\pi}}}\,.\EE
We first consider the case
\BE{ellannbd_e4} \big\|\ze|_{B_{2r,r}}\big\|_1\ge \frac{1}{25}\|\ze\|_1.\EE
Using polar coordinates, define $\ti\ze\!\in\!C^{\i}(B_{1,r};\R^k)$ by 
$$\ti\ze(\rho,\th)=\be(\rho)\ze(\rho,\th).$$
By H\"older's inequality and Lemma~\ref{loj_lmm},
$$\big\|\ze|_{B_{2r,r}}\big\|_1
\le \sqrt{3\pi}r\|\ti\ze\|_2 \le \sqrt{3\pi}r\|\d\ti\ze\|_1
\le \sqrt{3\pi}r\big(\|\d\ze\|_1+\|\be'\|_{C^0}\|\ze|_{B_{1,1/2}}\big\|_1\big).$$
Along with the assumptions \e_ref{ellannbd_e3} and \e_ref{ellannbd_e4},
this  implies the bound with
$$C=25\frac{\sqrt{3\pi} r}{1-24\sqrt{3\pi}\|\be\|_{C^0}r}
\le\frac{25}{48}\,.$$
Finally, suppose 
\BE{ellannbd_e6} \big\|\ze|_{B_{2r,r}}\big\|_1\le \frac{1}{25}\|\ze\|_1.\EE
Split the annulus $B_{1,r}$ into 3 wedges of equal area; 
split each wedge into a large convex outer portion and a small inner portion 
by drawing the line segment tangent to the circle of radius~$r$ 
and with the end points on the sides of the wedges $2r$
from the center as in Figure~\ref{ann_fig2}.
By~\e_ref{ellannbd_e6},
\BE{ellannbd_e7} A\equiv\big\|\ze|_{\cD_+}\big\|_1\ge \frac{8}{25}\|\ze\|_1\EE
for the outer piece $\cD_+$ of some wedge~$\cD$.
If
$$\bigg|\int_{\cD_+}\ze\bigg|\le \frac{3}{10}A\,,$$
then by Lemma~\ref{c0plane_lmm}, \e_ref{ellannbd_e3}, and H\"older's inequality,
\begin{equation*}\begin{split}
A&\le \frac{3}{10}A+ 
\frac{2\left(\frac{\sqrt{3}}{2}\right)^2}{\frac{\pi}{3}\left(1-\left(\frac{1}{96\sqrt{3\pi}}\right)^2\right)}
\int_{\cD_+}\!\int_{\cD_+}|\d_y\ze||y\!-\!x|^{-1}\d y\d x\\
&\le \frac{3}{10}A+ 
\frac{9}{2\pi}\cdot\frac{7\sqrt{2}}{9}\cdot2\pi\sqrt{3}
\int_{\cD}|\d_y\ze|\d y
\le \frac{3}{10}A+  7\sqrt{2\pi}\|\d\ze\|_2\,.
\end{split}\end{equation*}
Along with the assumption \e_ref{ellannbd_e7},
this  implies the bound with $C\!=\!125\sqrt{2\pi}/4$.
If 
$$\bigg|\int_{\cD_+}\ze\bigg|\ge \frac{3}{10}A\,,$$
then by \e_ref{ellannbd_e6}, \e_ref{ellannbd_e7}, and~\e_ref{ellannbd0_e1}, 
\begin{equation*}\begin{split}
A\le\big\|\xi|_{\cD}\big\|_1
&\le \big\|\ze|_{\cD}\big\|_1-\bigg|\int_{\cD}\ze\bigg|
+\int_0^{2\pi}\bigg|\int_r^1\ze(\rho,\th)\rho\d\rho\bigg|\d\th\\
& \le \bigg(A+\frac{1}{8}A\bigg)-\bigg(\frac{3}{10}A-\frac{1}{8}A\bigg)
+\sqrt{\frac{\pi}{2}}\|\d\ze\|_2
=\frac{19}{20}A+\sqrt{\frac{\pi}{2}}\|\d\ze\|_2\,.
\end{split}\end{equation*}
Along with the assumption \e_ref{ellannbd_e7},
this  implies the bound with $C\!=\!125\sqrt{2\pi}/4$.
Since $\be$ can be chosen so that $\|\be'\|_{C^0}\!<\!3$
(actually arbitrarily close to~2), 
comparing with \e_ref{ellannbd0_e2} for $R/r\!=\!144\sqrt{3\pi}$
we conclude that the claim holds with $C\!=\!125\sqrt{2\pi}/4$
for all~$r$.

\subsection{Bundle sections along smooth maps}
\label{sob_bounded_diff}

\noindent
Let $(M,g)$ be a Riemannian manifold and 
$(E,\lr{,},\na)$ a normed vector bundle with connection over~$M$.
If $u\!\in\!C^{\i}(\ti{B}_{R,r};M)$, $\xi\!\in\!\Ga(u;E)$, and $p\!\ge\!1$, 
let
$$\|\xi\|_p\equiv\bigg(\int_{\ti{B}_{R,r}}|\xi|^p\bigg)^{1/p}\,,\qquad
\|\xi\|_{p,1}\equiv\|\xi\|_p+\|\na^u\xi\|_p\,.$$

\begin{lmm}\label{pqplane_lmm}
If $(M,g)$ is a Riemannian manifold, 
$(E,\lr{,},\na)$ is a normed vector bundle with connection over~$M$,
and $p,q\!\ge\!1$ are such that $1\!-\!2/p\ge-2/q$, 
for every compact subset $K\!\subset\!M$ there exists 
$C_{K;p,q}\!\in\!\R^+$ with the following property.
If $R\!\in\!\R^+$, $r\!\in\![0,R)$, $u\!\in\!C^{\i}(\ti{B}_{R,r};M)$
is such that $\Im\,u\!\subset\!K$, 
and $\xi\!\in\!\Ga_c(u;E)$, then
$$\|\xi\|_q\le C_{K;p,q}R^{1-\frac{2}{p}+\frac{2}{q}}
\big(\|\na^u\xi\|_p+\|\xi\!\otimes\!\d u\|_p\big).$$
\end{lmm}

\noindent
{\it Proof:} Let $\exp\!:TM\!\lra\!M$ be an exponential-like map
and $\{U_i\!: i\!\in\![N]\}$  a finite open cover of $K$ 
such that the $g$-diameter of each set $U_i$ is at most $r_{\exp}^g(K)/2$.
Let $\{W_i\!: i\!\in\![N]\}$ be an open cover of $K$ such that 
$\ov{W}_i\!\subset\!U_i$.
Choose smooth functions $\eta_i\!: M\!\lra\![0,1]$ such that
$\eta_i\!=\!1$ on $W_i$ and $\eta_i\!=\!0$ outside of~$U_i$.
For each $i\!\in\![N]$, pick $x_i\!\in\!W_i$.
For each $z\!\in\!u^{-1}(U_i)\!\subset\ti{B}_{R,r}$, 
define $\ti{u}_i(z)\!\in\!T_{x_i}M$ and $\xi_i(z)\!\in\!E_{x_i}$ by
$$\exp_{x_i}\ti{u}_i(z)=u(z),~~|\ti{u}_i(z)|\!<\!r_{\exp}(x_i); 
\qquad \Pi_{\ti{u}_i(z)}\xi_i(z)=\xi(z).$$
For any $z\!\in\!B_{R,r}$, put $\ti{\xi}_i(z)\!=\!\eta_i(u(z))\xi_i(z)$.
Since $\ti{\xi}_i\!\in\!C^{\i}_c(\ti{B}_{R,r};E_{x_i})$,
by Corollary~\ref{pqplane_crl} there exists $C_{i;p,q}\!>\!0$ such that 
\BE{pqplane_e1}
\big\|\xi|_{u^{-1}(W_i)}\big\|_q= \big\|\ti\xi_i|_{u^{-1}(W_i)}\big\|_q 
\le \|\ti\xi_i\|_q \le C_{i;p,q}R^{1-\frac{2}{p}+\frac{2}{q}}\|\d\ti\xi_i\|_p\,.\EE
Since $\d\ti\xi_i=(\d\eta_i\!\circ\d u)\xi_i+(\eta\!\circ u)\d\xi_i$
on $u^{-1}(U_i)$ and vanishes outside of $u^{-1}(U_i)$, 
\BE{pqplane_e2}  \|\d\ti\xi_i\|_p
\le \big\|\d\xi_i|_{u^{-1}(U_i)}\big\|_p+C_i\|\xi_i\!\otimes\!\d u\|_p. \EE
On the other hand, by Corollary~\ref{ptder_crl}, if $u(z)\!\in\!U_i$
\BE{pqplane_e3} 
\Big|\na^u\xi|_z-\Pi_{\ti{u}_i(z)}\!\circ\!\d_z\xi_i\Big|
\le C_K|\d_zu||\xi(z)|.\EE
Combining equations~\e_ref{pqplane_e1}-\e_ref{pqplane_e3},
we obtain
$$\big\|\xi|_{u^{-1}(W_i)}\big\|_q\le \ti{C}_{i;p,q}R^{1-\frac{2}{p}+\frac{2}{q}}
\big(\|\xi\|_{p,1}+\|\xi\!\otimes\!\d u\|_p\big).$$
The claim follows by summing the last inequality over all $i$.

\begin{lmm}\label{c0plane_lmm2}
If $(M,g)$ is a Riemannian manifold, 
$(E,\lr{,},\na)$ is a normed vector bundle with connection over~$M$,
and $p\!>\!2$, 
for every compact subset $K\!\subset\!M$ there exists 
$C_{K;p}\!\in\!C^{\i}(\R^+;\R)$ with the following property.
If $R\!\in\!\R^+$, $r\!\in\![0,R/2]$,  $u\!\in\!C^{\i}(B_{R,r};M)$
is such that $\Im\,u\!\subset\!K$, and $\xi\!\in\!\Ga(u;E)$, then
$$\|\xi\|_{C^0}\le C_{K;p}(R)\big(\|\xi\|_{p,1}+\|\xi\!\otimes\!\d u\|_p\big).$$
\end{lmm}

\noindent
{\it Proof:} We continue with the setup in the proof of Lemma~\ref{pqplane_lmm}. 
By Corollary~\ref{c0plane_crl},  
$$\big\|\xi|_{u^{-1}(W_i)}\big\|_{C^0}\le 
\|\ti\xi_i\|_{C^0}\le C_{i;p}(R)\|\ti\xi_i\|_{p,1}
\le C_{i;p}(R)\big(\big\|\xi|_{u^{-1}(U_i)}\big\|_p+\|\d\ti\xi_i\|_p\big).$$
As above, we obtain
$$\|\d\ti\xi_i\|_p\le C_i\big(\|\na^u\xi\|_p+\|\xi\!\otimes\!\d u\|_p\big),$$
and the claim follows.

\begin{prp}\label{c0bound_prp}
If $(M,g)$ is a Riemannian manifold, 
$(E,\lr{,},\na)$ is a normed vector bundle with connection over~$M$,
and $p\!>\!2$, 
for every compact subset $K\!\subset\!M$ there exists 
$C_{K;p}\!\in\!C^{\i}(\R^+\!\times\!\R;\R)$ with the following property.
If $R\!\in\!\R^+$, $r\!\in\![0,R/2]$,  $u\!\in\!C^{\i}(B_{R,r};M)$
is such that $\Im\,u\!\subset\!K$, and $\xi\!\in\!\Ga_c(u;E)$, then
$$\|\xi\|_{C^0}\le C_{K;p}\big(R,\|\d u\|_p\big)\|\xi\|_{p,1}.$$
The same statement holds if $B_{R,r}$ is replaced by 
a fixed compact Riemann surface~$(\Si,g_{\Si})$.
\end{prp}

\noindent
{\it Proof:} By Lemma~\ref{c0plane_lmm2} applied with 
$\ti{p}=(p\!+\!2)/2$ and H\"older's inequality,
\BE{c0_plane_e1}
\|\xi\|_{C^0}\le C_{K;\ti{p}}(R)
\big(\|\xi\|_{\ti{p},1}+\|\xi\!\otimes\!\d u\|_{\ti{p}}\big)
\le \ti{C}_{K;\ti{p}}(R)\big(\|\xi\|_{p,1}+\|\d u\|_p\|\xi\|_{q_1}\big),\EE
where $q_1=p(p\!+\!2)/(p\!-\!2)$. 
If $q_1\!\le\!p$, then the proof is complete.
Otherwise, apply Lemma~\ref{pqplane_lmm}  with
$p_1=2q_1/(q_1\!+\!2)$ and H\"older's inequality:
\BE{c0_plane_e2}
\|\xi\|_{q_1}\le C_{K;p_1,q_1}(R)\big(\|\xi\|_{p_1,1}+\|\xi\!\otimes\!\d u\|_{p_1}\big)
\le C_{K;1}(R)\big(\|\xi\|_{p,1}+\|\d u\|_p\|\xi\|_{q_2}\big),\EE
where $q_2=pp_1/(p\!-\!p_1)$. 
If $q_2\!\le\!p$, then the claim follows
from equations~\e_ref{c0_plane_e1} and~\e_ref{c0_plane_e2}.
Otherwise, we can continue and construct sequences $\{p_i\},\{q_i\},\{C_{K;i}\}$
such that
\begin{gather}\label{c0_plane_e3}
p_i=\frac{2q_i}{q_i+2},\quad q_{i+1}=\frac{pp_i}{p-p_i};\\
\label{c0_plane_e4}
\|\xi\|_{q_i}\le C_{K;i}(R)\big(\|\xi\|_{p,1}+\|\d u\|_p\|\xi\|_{q_{i+1}}\big).
\end{gather}
The recursion~\e_ref{c0_plane_e3} implies that
$$q_{i+1}=\frac{2p}{2p+(p\!-\!2)q_i}q_i \quad\Lra\quad
\hbox{if~}q_i>0,\hbox{~then~~} 0<q_{i+1}<q_i.$$
Thus, if $q_i\!>\!2$ for all $i$, then 
the sequence $\{q_i\}$ must have a limit $q\!\ge\!2$ with 
$$q=\frac{2p}{2p+(p\!-\!2)q}q \quad\Lra\quad (p-2)q=0 \quad\Lra\quad q=0,$$
since $p\!>\!2$ by assumption.
Thus,  $q_N\!\le\!p$ for $N$ sufficiently large
and the first claim follows from \e_ref{c0_plane_e1}
and the equations~\e_ref{c0_plane_e4} with $i$ running from $1$ to~$N$,
where $N$ is the smallest integer such that~$q_{N+1}\!\le\!p$.
The second claim follows immediately from the first.

\subsection{Elliptic estimates}
\label{ellipticest_subs}

\noindent
If $A_1\!=\!B_{R_1,r_1}$ and 
$A_2\!=\!\bar{B}_{R_2,r_2}$
are two annuli in $\R^2$, we write $A_2\!\Subset_{\de}\! A_1$ if 
$R_1\!-\!R_2\!>\!\de$ and $r_2\!-\!r_1\!\ge\!\de$.

\begin{lmm}\label{elli_lmm1}
For any $\de\!>\!0$, $p\!\ge\!1$, and  open annulus~$A_1$, there exists 
$C_{\de,p}(A_1)\!>\!0$ such that for any annulus
$A_2\!\Subset_{\de}\! A_1$ and $\xi\!\in\!C^{\i}(A_1;\C^k)$,
$$\big\|\xi|_{A_2}\big\|_{p,1}\le C_{\de,p}(A_1)
\big(\|\dbar\xi\|_p+\|\d\xi\|_2+\|\xi\|_1\big),$$
where the norms are taken with respect to the standard metric on $\R^2$.
\end{lmm}

\noindent
{\it Proof:} We can assume that $A_2$ is the maximal annulus
such that $A_2\!\Subset_{\de}\! A_1$.
Let $\eta\!: A_1\!\lra\![0,1]$ be a compactly supported 
smooth function such that $\eta|_{A_2}\!=\!1$.
By the fundamental elliptic inequality for the $\dbar$-operator on $S^2$
\cite[Lemma~C.2.1]{MS},
\BE{elli1_e1}\begin{split}
\big\|\xi|_{A_2}\big\|_{p,1} &\le\|\eta\xi\|_{p,1} 
\le C_p(A_1)\big(\|\dbar(\eta\xi)\|_p\!+\!\|\eta\xi\|_p\big)\\
&\le C_p(A_1)
\big(\|\dbar\xi\|_p\!+\!\|(\d\eta)\xi\|_p\!+\!\|\eta\xi\|_p\big).
\end{split}\EE
By Corollary~\ref{pqplane_crl} with $(p,q)\!=\!(2,p)$ 
and $(p,q)\!=\!(1,2)$ and H\"older's inequality,
\BE{elli1_e2}\begin{split}
\|\eta\xi\|_p&\le C_p(A_1)\|\d(\eta\xi)\|_2
\le C_p(A_1)\big(\|\d\xi\|_2+\|(\d\eta)\xi\|_2\big)\\
&\le \ti{C}_p(A_1)\big(\|\d\xi\|_2+\|\d((\d\eta)\xi)\|_1\big)
\le \ti{C}_{p,\de}(A_1)\big(\|\d\xi\|_2+\|\d\xi\|_1+\|\xi\|_1\big)\\
&\le C_{\de,p}(A_1)\big(\|\d\xi\|_2+\|\xi\|_1\big).
\end{split}\EE
Similarly,
\BE{elli1_e3}
\|(\d\eta)\xi\|_p\le C_{\de,p}(A_1)\big(\|\d\xi\|_2+\|\xi\|_1\big).\EE
The claim follows by plugging \e_ref{elli1_e2} and \e_ref{elli1_e3}
into \e_ref{elli1_e1}.

\begin{crl}\label{elli_crl}
For any $\de\!>\!0$, $p\!\ge\!1$, and open annulus $A_1$, 
there exists $C_{\de,p}(A_1)\!>\!0$ such that for 
any annulus $A_2\!\Subset_{\de}\!A_1$, and 
$\xi\!\in\!C^{\i}(A_1;\C^n)$,
$$\|\d\xi|_{A_2}\|_p\le C_{\de,p}(A_1)
\big(\|\dbar\xi\|_p+\|\d\xi\|_2\big).$$
\end{crl}

\noindent
{\it Proof:} With $|A_1|$ denoting the area of $A_1$, let
$$\bar\xi=\frac{1}{|A_1|}\int_{A_1}\xi$$
be the average value of $\xi$.
By Lemma~\ref{elli_lmm1},
\BE{elli2_e1}\begin{split}
\|\d\xi|_{A_2}\|_p =\|\d(\xi\!-\!\bar\xi)|_{A_2}\|_p
&\le C_{\de,p}(A_1)\big(\|\dbar(\xi\!-\!\bar\xi)\|_p+
\|\d(\xi\!-\!\bar\xi)\|_2+\|\xi\!-\!\bar\xi\|_1\big)\\
&=C_{\de,p}(A_1)\big(\|\dbar\xi\|_p+\|\d\xi\|_2+\|\xi\!-\!\bar\xi\|_1\big).
\end{split}\EE
The claim follows by applying Corollary~\ref{ellannbd_crl}
with $\ze\!=\!\xi\!-\!\bar\xi$.

\begin{rmk}\label{ellicrl_rmk}
The case $r_1\!>\!0$ (which is the case needed for gluing pseudo-holomorphic maps
in symplectic topology) follows from Corollary~\ref{ellannbd_crl0};
Corollary~\ref{ellannbd_crl} can be used to obtain a sharper statement in this case
(that $C_{\de,p}(A_1)$ does not depend on~$r_1$).
The $r_1\!=\!0$ case requires only the first two steps
in the proof of Corollary~\ref{ellannbd_crl}.
\end{rmk}

\noindent
A \sf{smooth generalized CR-operator} in a smooth complex vector bundle
$(E,\na)$ with connection over an almost complex manifold $(M,J)$
is an operator of the form
$$D=\dbar_{\na}+A\!: \Ga(M;E)\lra \Ga(M;T^*M^{0,1}\!\otimes\!_{\C}E),$$
where 
$$\dbar_{\na}\xi=\frac{1}{2}\big(\na\xi+\fI\na_J\xi\big) 
\quad\forall\,\xi\!\in\!\Ga(M;TM),
\qquad A\in \Ga\big(M;\Hom(E;T^*M^{0,1}\!\otimes\!_{\C}E)\big).$$
If in addition $u\!:\Si\!\lra\!M$ is a smooth map from an almost complex manifold
$(\Si,\fj)$, the pull-back CR-operator is given~by
$$D_u=\dbar_{\na^u}+A\circ\partial u\!: \Ga(u;E)\lra \Ga^{0,1}(u;E).$$

\begin{prp}\label{elli_prp1}
If $(M,g)$ is a Riemannian manifold with an almost complex structure $J$, 
$(E,\lr{,},\na)$ is a normed complex vector bundle with connection over~$M$
and a smooth generalized CR-operator~$D$, and $p\!\ge\!1$, 
then for every compact subset $K\!\subset\!M$,  
$\de\!>\!0$, and open annulus $A_1\!\subset\!\R^2$, there exists 
$C_{K;\de,p}(A_1)\!\in\!\R^+$ with the following property.
If  $u\!\in\!C^{\i}(A_1;M)$ is such that $\Im\,u\!\subset\!K$,
$\xi\!\in\!\Ga(u;E)$, and 
$A_2\!\Subset_{\de}\!A_1$ is an annulus, then
$$\big\|\na^u\xi|_{A_2}\big\|_p\le 
C_{K;\de,p}(A_1) \big(\|D_u\xi\|_p+\|\na^u\xi\|_2+\|\xi\!\otimes\!\d u\|_p\big),$$
where the norms are taken with respect to 
the standard metric on $\R^2$.
\end{prp}

\noindent
{\it Proof:} We continue with the setup in the proof of Lemma~\ref{pqplane_lmm}. 
By Corollary~\ref{elli_crl},
\BE{elli3_e1}\begin{split}
\big\|\d\ti\xi_i|_{A_2}\big\|_p&
\le C_{i;\de,p}(A_1) \big(\|\dbar\ti\xi_i\|_p+\|\d\ti\xi_i\|_2\big)\\
&\le C_{i;\de,p}'(A_1)\big(\big\|\dbar\xi_i|_{u^{-1}(U_i)}\big\|_p
+\big\|\d\xi_i|_{u^{-1}(U_i)}\big\|_2+\|\xi\!\otimes\!\d u\|_p\big).
\end{split}\EE
Since $\na$ commutes with the complex structure in $E$ and 
$\ti\xi_i\!=\!\xi_i$ on $u^{-1}(W_i)$, 
it follows from \e_ref{pqplane_e3} and~\e_ref{elli3_e1} that
\BE{elli3_e3}\begin{split}
\big\|\na^u\xi|_{A_2\cap u^{-1}(W_i)}\big\|_p
&\le \big\|\d\ti\xi_i|_{A_2}\big\|_p+C_K\|\xi\!\otimes\!\d u\|_p\\
&\le \ti{C}_{i;\de,p}(A_1)
\big(\|\dbar_{\na^u}\xi\|_p+\|\na^u\xi\|_2+\|\xi\!\otimes\!\d u\|_p\big)\\
&\le \ti{C}_{i;\de,p}'(A_1)\big(\|D_u\xi\|_p+\|\na^u\xi\|_2+\|\xi\!\otimes\!\d u\|_p\big).
\end{split}\EE
The claim is obtained by summing the last equation over all~$i$.

\begin{lmm}\label{elli7}
If $(M,g)$ is a Riemannian manifold with an almost complex structure $J$, 
$(E,\lr{,},\na)$ is a normed complex vector bundle with connection over~$M$
and a smooth generalized CR-operator~$D$, and $p\!>\!2$, 
then for every compact subset $K\!\subset\!M$ and open ball $B\!\subset\!\R^2$,
there exists $C_{K;B,p}\!\in\!C^{\i}(\R;\R)$ 
with the following property.
If  $u\!\in\!C^{\i}(B;M)$ is such that $\Im\,u\!\subset\!K$ and
$\xi\!\in\!\Ga_c(u;E)$, then
$$\|\xi\|_{p,1}\le C_{K;B,p}(\|\d u\|_p)
\big(\|D_u\xi\|_p+\|\xi\|_p\big),$$
where the norms are taken with respect to the standard metric on $\R^2$.
\end{lmm}

\noindent
{\it Proof:} By an argument nearly identical to the proof of Proposition~\ref{elli_prp1},
$$\|\xi\|_{p',1}\le C_{K;p'}(B)\big(\|D_u\xi\|_{p'}+\|\xi\|_{p'}+
\|\xi\!\otimes\!\d u\|_{p'}\big)$$
for any $p'\!\ge\!1$.
On the other hand, by Proposition~\ref{c0bound_prp}, 
$$ \|\xi\|_{C^0} \le C_{K;B,\ti{p}}(\|\d u\|_{\ti{p}})\|\xi\|_{\ti{p},1},$$
where $\ti{p}\!=\!(p+2)/2$.
Proceeding as in the proof of Proposition~\ref{c0bound_prp}, we then obtain
\begin{equation*}\begin{split}
\|\xi\|_{p,1}&\le C_{K;B,p}(\|\d u\|_{\ti{p}})
\big(\|D_u\xi\|_p+\|\xi\|_p+\|\d u\|_p\|\xi\|_{\ti{p},1}\big),\\
\|\xi\|_{\ti{p},1}&\le C_{K;\ti{p}}(B)
\big(\|D_u\xi\|_p+\|\xi\|_p+\|\d u\|_p\|\xi\|_{q_1}\big),\\
\|\xi\|_{q_i}&\le C_{K;p_i,q_i}(B)\big(\|\xi\|_{p_i,1}+\|\xi\!\otimes\!\d u\|_{p_i}\big)\\
& \le C_{K;B,i}(\|\d u\|_p)\big(\|D_u\xi\|_p+\|\xi\|_p+\|\d u\|_p\|\xi\|_{q_{i+1}}\big);
\end{split}\end{equation*}
we stop the recursion at the same value of $i\!=\!N$ as in the proof of 
Proposition~\ref{c0bound_prp}.

\begin{prp}\label{elli_prp2}
If $(M,g)$ is a Riemannian manifold with an almost complex structure $J$, 
$(E,\lr{,},\na)$ is a normed complex vector bundle with connection over~$M$
and a smooth generalized CR-operator~$D$, and $p\!>\!2$, 
then for every compact subset $K\!\subset\!M$ and 
compact Riemann surface $(\Si,g_{\Si})$,
there exists $C_{K;\Si,p}\!\in\!C^{\i}(\R;\R)$ 
with the following property.
If  $u\!\in\!C^{\i}(\Si;M)$ is such that $\Im\,u\!\subset\!K$ and
$\xi\!\in\!\Ga(u;E)$, then  
$$\|\xi\|_{p,1}\le C_{K;\Si,p}\big(\|\d u\|_p)(\|D_u\xi\|_p+\|\xi\|_p\big).$$
\end{prp}

\noindent
{\it Proof:} This statement is immediate from Lemma~\ref{elli7}.\\

\vspace{.2in}

\noindent
{\it Department of Mathematics, SUNY Stony Brook, NY 11794-3651\\
azinger@math.sunysb.edu}\\

\end{document}